\documentclass[12pt]{amsart}

\usepackage{amssymb,amsmath,amsthm,amsfonts,amscd,bbm}
\usepackage{enumitem}
\usepackage{pdflscape}% rotating the page
\usepackage[paper=portrait,pagesize]{typearea}
\usepackage{caption}
\usepackage{bm}

\usepackage{geometry} % change margin

\usepackage{ifpdf}
\ifpdf
\usepackage[pdftex]{graphicx}
\else
\usepackage[dvips]{graphicx}
\fi
\usepackage[all]{xy}
\usepackage{hyperref}
\usepackage{xcolor}
\usepackage{graphicx}
\usepackage[outdir=./]{epstopdf}
\usepackage{enumitem}
\usepackage{tensor}
\usepackage{tikz-cd}
\usepackage{multicol}
\usepackage{mathtools}
\usepackage{tocvsec2}
\usepackage{bbm}
\usepackage{longtable}
\usepackage{fullpage}

\newcommand{\doi}[1]{\href{https://dx.doi.org/#1}{{\small DOI:#1}}}

\newcommand{\googlebooks}[1]{(preview at \href{https://books.google.com/books?id=#1}{google books})}

\newcommand{\numdam}[1]{}
\usepackage{mathrsfs}

\usepackage{mathrsfs}

\usepackage{mathrsfs}
\DeclareMathAlphabet{\mathpzc}{OT1}{pzc}{m}{it}

\def\semicolon{;}
\def\applytolist#1{
    \expandafter\def\csname multi#1\endcsname##1{
        \def\multiack{##1}\ifx\multiack\semicolon
            \def\next{\relax}
        \else
            \csname #1\endcsname{##1}
            \def\next{\csname multi#1\endcsname}
        \fi
        \next}
    \csname multi#1\endcsname}

\def\calc#1{\expandafter\def\csname c#1\endcsname{{\mathcal #1}}}
\applytolist{calc}QWERTYUIOPLKJHGFDSAZXCVBNM;
\def\bbc#1{\expandafter\def\csname bb#1\endcsname{{\mathbb #1}}}
\applytolist{bbc}QWERTYUIOPLKJHGFDSAZXCVBNM;
\def\bfc#1{\expandafter\def\csname bf#1\endcsname{{\mathbf #1}}}
\applytolist{bfc}QWERTYUIOPLKJHGFDSAZXCVBNM;
\def\sfc#1{\expandafter\def\csname s#1\endcsname{{\sf #1}}}
\applytolist{sfc}QWERTYUIOPLKJHGFDSAZXCVBNM;
\def\fc#1{\expandafter\def\csname f#1\endcsname{{\mathfrak #1}}}
\applytolist{fc}QWERTYUIOPLKJHGFDSAZXCVBNM;

\usepackage{tikz}
\usepackage{tikz-cd}

\makeatletter
\def\fixtikzforbreqn#1#2{%
  \protected\edef#1{\noexpand\ifmmode\mathchar\the\mathcode`#2 \noexpand\else#2\noexpand\fi}%
}
\fixtikzforbreqn\tikz@nonactivesemicolon;
\fixtikzforbreqn\tikz@nonactivecolon:
\fixtikzforbreqn\tikz@nonactivebar|
\fixtikzforbreqn\tikz@nonactiveexlmark!
\makeatother

\usepackage{breqn}   % Fix : | ; ! in tikzcd. And must put makealetter-other between tikzcd package and breqn package

\usetikzlibrary{arrows,backgrounds,patterns.meta}
\usetikzlibrary{positioning,shadings,cd}
\usetikzlibrary{shapes}
\usetikzlibrary{backgrounds}
\usetikzlibrary{decorations,decorations.pathreplacing,decorations.markings,decorations.pathmorphing}
\usetikzlibrary{fit,calc,through}
\usetikzlibrary{external}
\usetikzlibrary{arrows}
\tikzset{vertex/.style = {shape=circle,draw,fill=black,inner sep=0pt,minimum size=5pt}}
\tikzset{edge/.style = {->,> = latex', bend right}}
\tikzset{
	super thick/.style={line width=3pt}
}
\tikzset{
    quadruple/.style args={[#1] in [#2] in [#3] in [#4]}{
        #1,preaction={preaction={preaction={draw,#4},draw,#3}, draw,#2}
    }
}
\tikzstyle{shaded}=[fill=red!10!blue!20!gray!30!white]
\tikzstyle{unshaded}=[fill=white]
\tikzstyle{empty box}=[circle, draw, thick, fill=white, opaque, inner sep=2mm]
\tikzstyle{annular}=[scale=.7, inner sep=1mm, baseline]
\tikzstyle{rectangular}=[scale=.75, inner sep=1mm, baseline=-.1cm]
\tikzstyle{mid>}=[decoration={markings, mark=at position 0.5 with {\arrow{>}}}, postaction={decorate}]
\tikzstyle{mid<}=[decoration={markings, mark=at position 0.5 with {\arrow{<}}}, postaction={decorate}]
\tikzstyle{over}=[double, draw=white, super thick, double=]
\tikzstyle{snake}=[decorate, decoration={snake, segment length=1mm, amplitude=.3mm}]
\tikzstyle{saw}=[decorate, decoration={saw, segment length=.7mm, amplitude=.25mm}]

\tikzstyle{coupon}=[draw, very thick, rectangle, rounded corners=5pt]
\tikzset{Rightarrow/.style={double equal sign distance,>={Implies},->},
triplecd/.style={-,preaction={draw,Rightarrow}},
quadruplecd/.style={preaction={draw,Rightarrow,
shorten >=0pt
},
shorten >=1pt,
-,double,double
distance=0.2pt}}
\tikzset{
    tripleline/.style args={[#1] in [#2] in [#3]}{
        #1,preaction={preaction={draw,#3},draw,#2}
    }
}
\tikzstyle{triple}=[tripleline={[line width=.15mm,black] in
      [line width=.7mm,white] in
      [line width=1mm,black]}] 
\tikzset{
    quadrupleline/.style args={[#1] in [#2] in [#3] in [#4]}{
        #1,preaction={preaction={preaction={draw,#4},draw,#3}, draw,#2}
    }
}
\tikzstyle{quadruple}=[quadrupleline={[line width=.3mm,white] in
      [line width=.6mm,black] in
      [line width=1.2mm,white] in
      [line width=1.5mm,black]}]
 
\newcommand{\tikzmath}[2][]
     {\vcenter{\hbox{\begin{tikzpicture}[#1]#2
                     \end{tikzpicture}}}
     }

\newcommand{\roundNbox}[6]{
	\draw[rounded corners=5pt, very thick, #1] ($#2+(-#3,-#3)+(-#4,0)$) rectangle ($#2+(#3,#3)+(#5,0)$);
	\coordinate (ZZa) at ($#2+(-#4,0)$);
	\coordinate (ZZb) at ($#2+(#5,0)$);
	\node at ($1/2*(ZZa)+1/2*(ZZb)$) {#6};
}

\makeatletter
\newcommand{\xMapsto}[2][]{\ext@arrow 0599{\Mapstofill@}{#1}{#2}}
\def\Mapstofill@{\arrowfill@{\Mapstochar\Relbar}\Relbar\Rightarrow}
\makeatother

\newcommand{\AColor}{gray!30}
\newcommand{\BColor}{gray!55}
\newcommand{\CColor}{gray!80}

\theoremstyle{plain}
\newtheorem{thm}{Theorem}[section]
\newtheorem*{thm*}{Theorem}
\newtheorem{thmalpha}{Theorem}

\newtheorem*{cor*}{Corollary}

\newtheorem*{conj*}{Conjecture}
\newtheorem{lem}[thm]{Lemma}
\newtheorem*{lem*}{Lemma}

\newtheorem{prop}[thm]{Proposition}

\newtheorem*{quest*}{Question}
\newtheorem*{claim*}{Claim}

\theoremstyle{definition}

\newtheorem{ex}[thm]{Example}
\newtheorem*{ex*}{Example}
\newtheorem{defn}[thm]{Definition}

\newtheorem{construction}[thm]{Construction}

\newtheorem*{exs*}{Examples}
\newtheorem{sub-ex}[thm]{Sub-Example}
\newtheorem{counter-ex}[thm]{Counter-Example}

\newtheorem*{rem*}{Remark}
\newtheorem{remark}[thm]{Remark}

\usepackage{xcolor}
\definecolor{dark-red}{rgb}{0.7,0.25,0.25}
\definecolor{dark-blue}{rgb}{0.15,0.15,0.55}
\definecolor{medium-blue}{rgb}{0,0,.8}
\definecolor{DarkGreen}{RGB}{0,150,0}
\definecolor{rho}{named}{red}
\hypersetup{
   colorlinks, linkcolor={purple},
   citecolor={medium-blue}, urlcolor={medium-blue}
}

\newcommand{\id}{\operatorname{id}}

\newcommand{\Irr}{\operatorname{Irr}}

\newcommand{\op}{\operatorname{op}}

\newcommand{\tr}{\operatorname{tr}}
\newcommand{\Tr}{\operatorname{Tr}}

\newcommand{\End}{\operatorname{End}}

\newcommand{\Aut}{\operatorname{Aut}}

\newcommand{\fgpBim}{{\sf{Bim}_{\sf fgp}}}
\newcommand{\TLJ}{{\sf TLJ}}
\newcommand{\QSet}{{\sf QSet}}

\newcommand{\Rep}{{\mathsf{Rep}}}

\newcommand{\Hilb}{\mathsf{Hilb}}
\newcommand{\fdHilb}{{\sf Hilb_{fd}}}

\newcommand{\rCorr}{{\mathsf{C^*Alg}}}
\newcommand{\vNAlg}{{\mathsf{vNAlg}}}

\newcommand{\Graph}{\mathsf{Graph}}

\newcommand{\For}{\mathsf{For}}

\def\altdb{\vadjust{\vbox to 0pt{\vss\hbox{\kern \hsize
\quad{\dbend}}\kern\baselineskip\kern-10pt}}}

\newcommand{\rip}[3]{\langle #1 \ |\ #2\rangle_{#3}}
\newcommand{\lip}[3]{{_{#1}}\langle #2, #3\rangle}
\newcommand{\ev}[1]{\mathsf{ev}_{#1}}
\newcommand{\coev}[1]{\mathsf{coev}_{#1}}

\newcommand{\noshow}[1]{}
\renewcommand{\MR}[1]{}

\title{ Quantum graphs, subfactors and tensor categories I}
\date{}
\author{Michael Brannan and Roberto Hern\'{a}ndez Palomares}

\begin{document}
\begin{abstract}
    We develop an equivariant theory of graphs with respect to quantum symmetries and present a detailed exposition of various examples. 
    We portray unitary tensor categories as a unifying framework encompassing all finite classical simple graphs, (quantum) Cayley graphs of finite (quantum) groupoids, and all finite-dimensional quantum graphs. 
    We model a quantum set by a finite-index inclusion of C*-algebras and use the quantum Fourier transform to obtain all possible adjacency operators.
    In particular, we show every finite-index subfactor can be regarded as a complete quantum graph and describe how to find all its subgraphs. 
    As applications, we prove a version of Frucht's Theorem for finite quantum groupoids, and introduce a version of path spaces for quantum graphs. 
\end{abstract}

\maketitle
\section{Introduction}
Quantum graphs are nowadays ubiquitous objects of study at the interface of quantum information theory, operator algebras and category theory. 
They arise in very many flavors and contexts, e.g. as confusability graphs of noisy quantum channels in zero-error communication problems \cite[\S7]{MR3549479}, as natural models for completeness problems in quantum complexity theory \cite{Culf-Mehta}, and in the study of non-local games \cite{MR3849575, MR3907958, MR4091496, MR4505911, MR4752739, MR4507619, 2024arXiv240815444G}. 
Quantum graphs also appear as the Cayley graphs of quantum groups \cite{MR2130588, 2023arXiv230615315W}, as well as examples of quantum relations on von Neumann algebras \cite{MR2908249, MR4302212}, and have recently been studied  as objects of interest in their own right.

Quantum graphs have appeared recently in \cite{MR3849575, MR3907958, MR4481115, MR4514486, MR4555986, 2023arXiv230615315W,  MR4706978}, as formulated in more operator algebraic and categorical terms. 
Inspired by Gelfand duality and its noncommutative generalizations, Musto, Reutter and Verdon \cite{MR3849575, MR3907958} described a theory of  finite quantum graphs in terms of certain Frobenius algebras in the category of finite dimensional Hilbert spaces, encompassing all quantum graphs over finite dimensional C*-algebras (equipped with certain canonical tracial positive form -- a condition which was later relaxed in the works \cite{MR4091496, MR4481115, MR4514486, MR4706978}).  
In this work, we further broaden the above viewpoint,  laying out a theory of quantum graphs internal to $\cC,$ a unitary tensor category (UTC). 
Doing so allows one to frame abundant families of examples such as classical graphs and Cayley graphs of finite groupoids and their quantum analogues, all under the same umbrella.   Since UTCs $\cC$ arise most naturally as (quantum) symmetries (e.g., representation categories of quantum groups, standard invariants of subfactors, and so on), our approach can be viewed as a $\cC$-equivariant theory of (finite, simple, directed, quantum) graphs.

We define a $\cC$-equivariant graph in terms of a Q-system in $\cC$, together with a chosen convolution idempotent $\hat{T}\in \End_\cC(Q)$.
Whenever $\cC$ is not explicitly given or is not inferred from context, we refer to $\cC$-equivariant graphs simply as categorified graphs. 
A Q-system in $\cC$ can be viewed as a categorification of a finite dimensional C*-algebra equipped with a special state  (\emph{aka} $\delta$-form \cite{MR1889999} for certain diagrammatic categories).  In our framework, the Q-system $Q$ plays the role of a noncommutative space of vertices (a {\it finite quantum set}), while the morphism $\hat{T}\in \End_\cC(Q)$ quantizes the  adjacency data. In case $\cG=(V,E)$ is a classical graph,  $Q = C(V)$, and $\hat T$ is just the adjacency matrix $\hat{T}\in \mathsf{Mat}_{V\times V}(\{0,1\}),$ where $\hat{T}(u,v)=1$ if and only if $(u,v)\in E.$ All the graphs considered in this paper are therefore simple; this is, having no ``multiple edges''. 

The UTCs that one encounters in nature are often \emph{concrete}, as for any $\cC$ there is a unital simple C*-algebra $A$ with $\cC\subset \fgpBim(A)$ as a full $\otimes$-subcategory. 
Here, $\fgpBim(A)$ is the UTC of finitely generated projective $A$-$A$ bimodules.
As a consequence, we concretely model a Q-system by a finite Jones/Watatani-index inclusion $A\subset B$ with a faithful expectation $E:B\twoheadrightarrow A$ \cite{MR696688,MR996807}.
A concrete $\cC$-equivariant graph then consists of an \emph{adjacency operator} $\hat{T}\in\End_{A-A}( B),$ the finite-dimensional C*-algebra of adjointable $A$-$A$ bimodular endomorphisms of $B.$
Here, $\hat T$ is a \emph{Schur idempotent}; i.e. $\hat T\star \hat T =  \hat T,$ where $\star$ denotes the Schur product (Equation \ref{eqn:SchurProduct}),  generalizing the convolution from groups, and also the Hadamard (i.e., entrywise) product of matrices. 
In this sense, all the concrete graphs we consider are thus ``finite''.

This viewpoint allows one  to incorporate tools from the theory of subfactors such as their standard invariants, the quantum Fourier transform \cite{Ocn90, MR4236188}, and reconstruction techniques extending the familiar Tannakian theory from groups.  
It is this tool-set we exploit throughout this work, allowing us to view inclusions as quantum sets, and moreover to parameterize the space of graphs that they can support.  
We make this precise in the following theorem, summarizing Sections \ref{sec:SubfacIndexStdInv}, \ref{sec:QFourierTransf}, and \ref{sec:Examples}: 
\begin{thmalpha}[{\bf Every categorified graph is concrete}]
    For any UTC $\cC$, every Q-system $Q\in \cC$ corresponds to some finite-index inclusion $A\subset B,$
    and each quantum adjacency operator on $Q$ to some  $\star$-idempotent $\hat{T}\in \End_{A-A}(B)$.
    
    Furthermore, the quantum Fourier transform 
    $$
        \cF:\End_{B-B}(B_1)\to \End_{A-A}(B),
    $$
    is a linear isomorphism exchanging composition $\circ$ for the $\star$-product, where $B_1\cong \End_{\bbC-A}(B)$ is the Jones basic construction $A\subset B\subset B_1.$
    The map $\cF$ restricts to a bijection between the $\circ$-idempotents and the $\star$-idempotents; and therefore every adjacency operator on $B$ is the Fourier transform of a $\circ$-idempotent on $B_1.$ 
\end{thmalpha}

\noindent The proof of this is essentially that every UTC can always be made concrete, constructing $A$ solely out of $\cC$. (cf \cite{MR4139893, MR4419534})
In this framework, the Jones projection $e=\cF(\id_{B_1})\in \End_{A-A}(B)$ whose range is $A$, is the adjecency operator of {\bf the complete graph on the inclusion $A \subset B$}. 
We refer to $A$ as the algebra of scalars, which in the classical case collapses to $A=\bbC.$ 

\medskip

At the very least, our setting should expressely recollect all classical graphs, and the finite-dimensional noncommutative graphs well-studied in  the literature.  We address this by highlighting some important families of  examples.
\begin{exs*}[{\bf All simple Finite-dimensional graphs}]\ \\ 
    $\circ$ \underline{All classical graphs} (Example \ref{ex:ClassicalGraphs}):  
    The scalars are $A=\bbC,$ and the vertices   correspond to the minimal projections in $B={\bbC}^N.$ 
    The conditional expectation is the  normalized trace $\tr$  on $\mathsf{Mat}_N(\bbC)$, where we view $B \subset \mathsf{Mat}_N(\bbC)$ as the diagonal subalgebra.  
    The Jones basic construction is $B_1 = \mathsf{Mat}_N(\bbC)$ considered with the same trace, and the first relative commutant $\End(\bbC^N)\cong \mathsf{Mat}_N(\bbC)$ is spanned by $\mathsf{Mat}_N(\{0,1\})$, the Schur-idempotent matrices--recovering all finite simple directed classical graphs. 

    \noindent$\circ$ \underline{All the ``usual'' matrix quantum graphs} (Example \ref{ex:UsualQuGraphs}): Starting from $A=\bbC I_N\subset \mathsf{Mat}_N(\bbC)=B$ equipped with $E=\omega,$ a faithful state, our methods recover all Schur idempotents on this quantum set. 
    Furthermore, if $N=2,$ and $\omega=\tr,$ this recovers the quantum graphs constructed and classified by Matsuda and Gromada \cite{MR4481115, MR4514486}.
\end{exs*}

We shall now exhibit \emph{increasingly quantum} families of examples of $\cC$- equivariant graphs:  

\begin{exs*}[{\bf Cayley graphs, spin models and subfactors}]\ \\
    $\bullet$ \underline{Finite groups} (Examples \ref{ex:GroupHyperfiniteCommutants}, \ref{ex:GammaGraphs} and \ref{ex:GroupsAndCayleyGraphs}):
    Given a finite group $\Gamma$, there is a hyperfinite subfactor 
    $$
        N:=\cR\subset \cR\rtimes\Gamma:=M
    $$
    of index $|\Gamma|,$ which is generated by the canonical \emph{commuting square}
    \begin{equation*}
    \begin{tikzcd}[column sep=0em, row sep=-.1em]
     \End_{N-N}(M)\cong \ell^\infty(\Gamma) & \subset   & \End_{N-N}(M_1)\cong \mathsf{Mat}_\Gamma(\bbC)\\
    \cup & {} & \cup \\
     \End_{M-M}(M)\cong \bbC & \subset &  \End_{M-M}(M_1)\cong \bbC[\Gamma] 
     \end{tikzcd}.
    \end{equation*}
    In this setting, the subfactor Fourier transform reduces to the non-abelian Fourier transform $\mathcal F:  (\bbC[\Gamma], \star) \to (\ell^\infty(\Gamma), \cdot)$.  Here, $\star$ denotes the convolution product on the group algebra, and $\cdot$ denotes the pointwise product of functions on $\Gamma$. Categorified graphs on $N \subset M$ then correspond to idempotents in $\bbC[\Gamma]$, which can be interpreted as {\it quantum Cayley graphs} on the quantum group $\hat \Gamma$ in the sense of \cite{2023arXiv230615315W}. (A slight difference is that our graphs are allowed to be \emph{degenerate}, meaning that $S$ need not generate $\Gamma$.)
    The usual notion of Cayley graphs on $\Gamma$ dually correspond to projections in $\ell^\infty(\Gamma)$.  These correspond to quantum graphs on the {\it dual subfactor} $M \subset M_1$ arising from the basic construction.

$\bullet$ \underline{{Spin models}}:
A square of algebras of the form 
\begin{equation*}
    \begin{tikzcd}[column sep=0em, row sep=-.1em]
     u\bbC^N u^* & \subset   & \ \mathsf{Mat}_N(\bbC)\\
    \cup & {} & \cup \\
     \bbC & \subset &   \bbC^N 
     \end{tikzcd}
    \end{equation*}
    satisfies the commuting square condition if and only if $u\in \End(\bbC^N)$ is a \emph{Hadamard matrix} \cite[\S5.2.2]{MR1473221}. That is, $u$ is a unitary whose entries all have the same modulus, which equivalently guarantees orthogonality between the MASAs $\bbC^N$ and $u\bbC^N u^*$ inside of $\mathsf{Mat}_N(\bbC).$ 
    Iterating this square gives a finite-index hyperfinite subfactor $N\subset M$ with {\it both} $N'\cap M_1$ and $M'\cap M_2$ abelian.    
    Using the Fourier transform one can then obtain all $\star$-idempotents on the quantum set ${}_NM_N$ from the projections in $M'\cap M_2.$ 
    The minimal projections in $M'\cap M_2$ can be conveniently described combinatorially in terms of a dual pair of association schemes associated to the Hadamard matrix $u$ \cite{MR1635553}.  We note that association schemes can be  interpreted as generalizations of the familiar Pontryagin duality from finite abelian groups. 
    We will further analyze these examples in the forthcoming articles \cite{BHPII, BGHP}.

$\bullet$ \underline{{Finite-index subfactors}} (Example \ref{ex:BiprojSubalgs}, Construction \ref{const:CompleteTrivQuGraph}, and Example \ref{Ex:HypGraphsLeq4}):  The broadest type of graphs equivariant with respect to a (concrete) UTC arise from finite-index subfactors. 
Associated to such a subfactor $N\subset M$ is its standard invariant, given by the lattice of higher relative commutants $\{N'\cap M_n\}_{n\geq 1}$ and $\{M'\cap M_n\}_{n\geq 1},$ where the $M_n$ are the steps of the Jones tower (c.f. Section \ref{sec:SubfacIndexStdInv}). 
Then all possible adjacency operators on the quantum set ${}_M{M_1}_M$ come from the $\circ$-idempotents in $\End_{N-N}(M)\cong N'\cap M_1$ through the quantum Fourier transform. Similarly, all adjacency operators on ${}_N{M}_N$ come from $\circ$-idempotents in $\End_{M-M}(M_1)\cong M'\cap M_2.$
 Since these algebras are finite-dimensional C*-algebras, it is in principal possible to obtain on-demand data about it's internal equivariant graphs (cf \cite{MR2718953}).
In particular, the Jones projection $e$ gives the complete quantum graph over $_NM_N$. 

In Example \ref{Ex:HypGraphsLeq4} we exhibit all quantum graphs supported on the small hyperfinite ADE-type subfactors. 
Addressing a question of Matsuda \cite{Junichiro}, we explain how the  $A_n$ subfactors constructed by Jones \cite{MR696688} give examples of \emph{regular} tracial quantum graphs of non-integer degree given by the square root of the index, which is impossible in tracial finite-dimensional quantum graphs \cite{Junichiro}.  
\end{exs*}
Similar constructions can be carried out in the representation category of a compact quantum group $\bbG$, where we can make use of the Fourier transform to define quantum Cayley graphs on the Q-systems in $\cC=\Rep(\bbG)$ (cf Example \ref{ex:Q-SysCompactQuantumGroups}).
This overlaps with the quantum Cayley graphs in \cite{2023arXiv230615315W}, defined over locally compact quantum groups.

\medskip
    
We now recall {\bf Frucht's Theorem} \cite{MR1557026}, which states that every finite group arises as the automorphism group of some finite simple undirected graph. 
While it is known that not every finite quantum group (= finite dimensional Hopf C$^\ast$-algebra) arises as the quantum automorphisms of some graph -- a simple counterexample is $\widehat{S_3}$, the Pontryagin dual of the permutation group $S_3$ \cite{MR4462380} -- one could ask if $\cC$-equivariant graphs are abundant enough to fill this gap. Actually the most natural conjecture here is that every finite quantum group arises as the quantum automorphism group of a {\it finite-dimensional quantum graph} (i.e., it suffices to take  $\cC  = \fdHilb$).  While this version of a quantum Frucht theorem seems out of reach at this time, we prove a variant of Frucht's Theorem for certain finite quantum groupoids, relying on a 2-category of $\cC$-equivariant graphs analogous to that of \cite{MR3849575}, whose morphisms are \emph{quantum graph homomorphisms}---``entanglement assisted'' compatible functions---, and their  intertwiners. 
\begin{thmalpha}[{\bf a quantum variant of Frucht's theorem} (Theorem \ref{thm:QuFrucht}) ]\label{thmalpha:Fruchts}
    Any quasitriangular finite quantum groupoid arises as the ``crossing quantum automorphisms'' of some $\cC$-equivariant graph. 
\end{thmalpha}
The first step in the proof of Theorem \ref{thm:QuFrucht} relies on the reconstruction results due to Nikshych and Vainerman, roughly, showing that finite quantum groupoids correspond to hyperfinite subfactors of finite index and finite depth (i.e., the associated category has finitely many simple objects) \cite{MR1800792}. 
The braiding on the associated subfactor afforded by quasitriangularity then gives the alluded crossing quantum automorphisms of the underlying quantum set. 
The last step is to prove that forgetting the braiding and retaining the entanglement bimodule is a monoidal equivalence between the category of crossing quantum automorphisms and the standard invariant for the subfactor.

\medskip
The language of subfactors adapts naturally to talk about paths in categorified graphs. 
For any $\cC$-equivariant graph $\cG$ we can make sense of the {\bf quantum edge correspondence} from \cite{MR4555986}, as the range of the \emph{edge projector} $E_\cG$ (equation \eqref{eqn:EdgeProjector}). 
In fact, using the graphical calculus in the C* 2-category of C*-correspondences, $E_\cG$ is quite transparently related to the $\circ$-idempotent $T$ whose Fourier transform $\hat{T}$ gives the adjacency operator defining $\cG$.
Explicitly, if we model a categorified graph as $\cG = (A\subset B,\  \hat{T}\in \End_{A-A}(B)),$ then $T\circ T = T\in \End_{B-B}(B_1),$ and so the edge correspondence is a finitely generated $B$-$B$ sub-bimodule $T[B_1]\subset B_1$ of the basic construction. 
Importantly, the relative tensor product $B\boxtimes_A B$ is isomorphic to $B_1$ as a C*-algebra.
And so, if $\hat{T}=e$ is the Jones projection yielding the complete graph, then the edges become $B\boxtimes_A B$. 
At the other extreme, when  $\hat{T}= \id_B$ (yielding the the trivial graph), then $T$ is the Jones projection from $B_1$ onto $B$, and so the edge correspondence trivializes to $B$, as one would expect for a trivial graph (with only self loops). 

As we continue along the Jones Tower $A\subset B\subset B_1\subset B_2\subset \cdots$, the construction of the edge correspondence suggests an organic definition of the higher {\bf quantum path spaces} $E^{(n)}_\cG$ as the range of the idempotent $T^{\boxtimes_B n}$ for $n\geq 1,$ as $B$-$B$ sub-bimodules of $B_n\cong B^{\boxtimes_A n}.$   
With this framework, we can recover some examples of interest developed in \cite[\S4]{MR4555986}: 
If $\cG$ is the complete graph, then $E^{(n)}_\cG= B_n$, and taking an inductive limit, $E^{(\bbN)}_\cG = B_\infty.$ 
In forthcoming articles we will study combinatorial and dynamical aspects of this construction, and its relation with subfactors. 

We mention that the C* 2-category of $\cC$-equivariant graphs and quantum graph homomorphisms as above is an interesting object of study, beyond $\cC=\fdHilb$ introduced in \cite{MR3907958}. It seems possible to use the $K$-theory of C*-algebras to introduce meaningful invariants for categorified graphs in the spirit of \cite{MR4717816}. And conversely, that non-local games such as the quantum graph homomorphism game  suggest the notion of \emph{values for games} for categorified graphs with coefficients in $A$. 
In \cite{BHPII, BGHP}, 
we shall explore these problems, including the connections with Weaver's quantum relations \cite{MR2908249} and the Hadamard graphs and matrices from spin models.

\tableofcontents

\subsection*{Acknowledgements}
The authors are very grateful to 
Moritz Weber and 2022 Focus semester in quantum information at Universit\"at des Saarlandes where this project was started.
This paper benefited from fruitful conversations with Daniel Gromada, Dietmar Bisch, Junichiro Matsuda, Brent Nelson, and  Corey Jones.
MB and RHP were partially supported by an NSERC Discovery Grant and NSF Grant DMS-2001163. RHP was partially supported by an AMS-Simons Travel Grant 2022.

\section{Preliminaries}\label{sec:Preliminaries}

\subsection{Categorical preliminaries and examples}
\ \\
In this section we briefly summarize the necessary background on unitary tensor categories and C* 2-categories, mainly to establish the notation and for the sake of self-containment. We assume nevertheless that the reader has some familiarity with (tensor) categories. For more details, the interested reader can consult \cite{MR4419534} and references therein, and \cite{MR3242743} for fundamentals of tensor categories. Regardless, we shall describe concrete illustrative examples of all these notions using C*-algebras.  

\begin{defn}
A \emph{$\rm C^*$-category} $\cC$ is a $\bbC$-linear category such that:
\begin{enumerate}
\item for all objects $a,b\in\cC$, there is a conjugate-linear involution
$$\dag:\cC(a\to b)\to \cC(b\to a)$$ satisfying $(f\circ g)^\dag=g^\dag\circ f^\dag.$ 
\item There is a Banach norm on morphisms such that $\|f\|^2 = \|f\circ f^\dag\|= \|f^\dag\circ f\|,$ for all $f\in \cC(a\to b)$.
\item For each $f\in \cC(a\to b)$, there exists $g\in \cC(a\to a)$ such that $f^\dag\circ f = g^\dag\circ g$.
\item
$\cC$ is Cauchy complete: there exists unitary direct sums, and every projection splits (see \cite[Assumption 2.7]{MR4419534}).
\end{enumerate}
\end{defn}

{\bf 2-categories.} All of our \emph{2-categories} are weak 2-categories, also known as a bicategories. We refer the reader to \cite{MR4261588} for background on 2-categories.

Given a 2-category $\cC$, its objects are denoted by lower case letters $a,b,c,\dots$, its 1-morphisms $a\to b$ are denoted by ${}_aX_b$, and 2-morphisms are by $f,g,h,\dots$.
We write 1-composition $\boxtimes$ from \emph{left to right}, e.g.~${}_aX\boxtimes_b Y_c$,
and 2-composition $\circ$ from \emph{right to left}.
This notation is consistent with the tensor product of bimodules and composition of intertwiners in our examples. 

\begin{defn}
    A \emph{$\rm C^*$ 2-category} is a 2-category such that every every $\hom$ 1-category is equipped with the structure of a $\rm C^*$-category, and the horizontal composition of 2-morphisms is compatible with the $\dag$-structures.
\end{defn}

We will often use the graphical calculus of string diagrams for 2-categories \cite[\S8.1.2]{MR3971584}. In a string diagram, shaded regions correspond to objects, 1-morphisms correspond to strings, and 2-morphisms correspond to coupons. For example, 

\begin{align*}
f: {}_aX\boxtimes_b Z_c \Rightarrow {}_aW_c
\qquad\rightsquigarrow\qquad
\tikzmath{
\begin{scope}
\clip[rounded corners=5pt] (-.7,-.7) rectangle (.7,.7);
\fill[\AColor] (-.7,-.7) -- (-.2,-.7) -- (-.2,0) -- (0,0) -- (0,.7) -- (-.7,.7);
\fill[\CColor] (.7,-.7) -- (.2,-.7) -- (.2,0) -- (0,0) -- (0,.7) -- (.7,.7);
\fill[\BColor] (-.2,-.7) rectangle (.2,0);
\end{scope}
\draw (0,0) -- (0,.7) node[above]{$\scriptstyle W$};
\draw (-.2,-.7) node[below]{$\scriptstyle X$} -- (-.2,0);
\draw (.2,-.7) node[below]{$\scriptstyle Z$} -- (.2,0);
\roundNbox{fill=white}{(0,0)}{.3}{.1}{.1}{$f$}
\node at (-.55,0) {$\scriptstyle a$};
\node at (.55,0) {$\scriptstyle c$};
\node at (0,-.55) {$\scriptstyle b$};
}    
\end{align*}
\noindent Our convention is to read the horizontal 1-composition from \emph{left to right}, and the vertical 2-composition from \emph{bottom to top}. We will also drop the sub-indices $a,b,c,...$ whenever no ambiguity arises.
The tensor product, composition and adjoint also admit graphical representations as follows: 
\begin{align*}
(f^\dag \boxtimes h)\circ g: \qquad
\tikzmath{
\begin{scope}
\clip[rounded corners=5pt] (-.7,-1.7) rectangle (1.7,.7);
\fill[\AColor] (-.7,.7) -- (-.2,.7) -- (-.2,0) -- (0,0) -- (0,-1.7) -- (-.7,-1.7);
\fill[\CColor] (2,.7) -- (.2,.7) -- (.2,0) -- (0,0) -- (0,.-7) -- (2,-.7);
\fill[\BColor] (-.2,.7) rectangle (.2,0);
\end{scope}
\draw (0,0) -- (0,-.7);
\draw (-.2,.7) node[above]{$\scriptstyle X$} -- (-.2,0);
\draw (.2,.7) node[above]{$\scriptstyle Z$} -- (.2,0);
\draw (1.1,.7) node[above]{$\scriptstyle Y$} -- (1.1,0);
\draw(1.1,0) -- (1.1, -.8);
\roundNbox{fill=white}{(0,0)}{.3}{.1}{.1}{$f^\dag$}
\roundNbox{fill=white}{(1.1,0)}{.3}{.1}{.1}{$h$}
\roundNbox{fill=white}{(.9,-1)}{.3}{.8}{.1}{$g$}
\draw (0,-1.7) node[below]{$\scriptstyle Z'$} -- (0,-1.3);
\draw (1.1,-1.7) node[below]{$\scriptstyle Y'$} -- (1.1,-1.3);
\node at (-.55,-.5) {$\scriptstyle a$};
\node at (1.55,-.5) {$\scriptstyle c$};
\node at (0,.55) {$\scriptstyle b$};
}.
\end{align*}
Here, $g: {}_aZ'\boxtimes_c Y'_c \Rightarrow {}_aW\boxtimes_c W'_c$, and $h:{}_cW'_{c}\Rightarrow {}_cY_c$.

The following is the prototypical example of a concrete C* 2-category, and we now lay a brief summary of its fundamental relevant features. An interested reader can consult further details in \cite[\S 1.1, 1.3]{2023arXiv230505072H}.
\begin{ex}[{\bf The C* 2-category of unital C*-algebras}]\label{ex:rCorr}
    The $\rm C^*$-2-category $\rCorr$ of right correspondences has objects unital $\rm C^*$-algebras $A,B,C,\hdots$, 1-morphisms ${}_AX_B, {}_AY_B,$ ${}_BZ_C$ are \emph{right correspondences}, and 2-morphisms ${}_AX_B \Rightarrow {}_AY_B$ are adjointable $A-B$ bimodular maps.

    For unital $C^*$-algebras $A,B$, {\bf a right C* correspondence}  ${}_AX_B\in\rCorr(A\to B)$ is a $\bbC$-vector space $X$ equipped with commuting left $A$- and right $B$-actions (written $a\rhd \xi$ and $\xi \lhd b$ for $a \in A, b \in B, \xi \in X$), and a right $B$-valued inner product $\langle \,\cdot\,|\,\cdot\,\rangle_B : \overline{X}\times X \to B$ which is:
    \begin{itemize}
    \item 
    right $B$-linear: $\langle \xi_1| \xi_2 \lhd b + \xi_3\rangle_B = \langle \xi_1| \xi_2 \rangle_B b + \langle \xi_1|\xi_3\rangle_B$,
    \item
    left conjugate $B$-linear: 
    $\langle \xi_1\lhd b + \xi_2 | \xi_3\rangle_B = b^*\langle \xi_1| \xi_2 \rangle_B + \langle \xi_2|\xi_3\rangle_B$
    \item
    positive: $\langle \xi|\xi\rangle_B \geq 0$ in $B$, and
    \item
    definite: $\langle \xi|\xi\rangle_B = 0$ implies $\xi=0$.
    \end{itemize}
    Observe that $\langle \eta|\xi\rangle_B^* = \langle \xi|\eta\rangle_B$ by the polarization identity, and we have the \emph{Cauchy-Schwarz inequality}:
    $$
    \langle \eta|\xi\rangle_B \langle\xi|\eta \rangle_B
    \leq
    \|\langle \xi|\xi\rangle_B\|_B \cdot \langle \eta|\eta\rangle_B
    \qquad\qquad
    \forall\, \eta,\xi\in X.
    $$
    This identity implies that
    $$
    \|\xi\|_X^2:=\|\langle \xi|\xi\rangle_B\|_B
    $$
    gives a well-defined norm on $X$. 
    We require that $X$ is complete with respect to the metric induced by this norm. 

    Moreover, we require the left $A$-action on $X$ is by \emph{adjointable operators}. 
    A right $B$-linear map $T:X_B \to Y_B$ between right $B$-modules is \emph{adjointable} if there is a right $B$-linear map $T^\dag:Y_B \to X_B$ such that 
    $$
    \langle \eta| T\xi\rangle_B = \langle T^\dag \eta|\xi\rangle_B
    \qquad\qquad
    \forall\, \xi\in X,\,\,
    \forall\, \eta\in Y.
    $$
    Observe that adjointable maps are necessarily bounded by the Closed Graph Theorem.

    The composition of 1-morphisms is given by the \emph{relative tensor product} $X\boxtimes_B Z$, which is defined as follows.
    First, we take the algebraic tensor product
    $X\otimes_\bbC Z$
    and consider the right $C$-valued sesquilinear form
    $$
    \langle \xi_1 \otimes \eta_1 | \xi_2\otimes \eta_2\rangle_C
    :=
    \langle  \eta_1 | \langle \xi_1|\xi_2\rangle^X_B \rhd \eta_2\rangle^Y_C.
    $$
    We let $N$ denote the left kernel of this form, i.e.,
    $$
    N=\{\zeta \in X\otimes_\bbC Z|\ \langle \zeta|\zeta\rangle_C = 0\},
    $$
    which is an algebraic $A-C$ sub-bimodule of $X\otimes_\bbC Z$,
    and note
    $$
    \operatorname{span}_{\bbC}\{(\xi\lhd b)\otimes_\bbC \eta - \xi\otimes_\bbC(b\rhd \eta)\}
    \subset N.
    $$
    Observe that $\langle \,\cdot\,|\,\cdot\,\rangle_C$ descends to a well-defined $C$-valued inner product on $(X\otimes_\bbC Y)/N$.
    We define $X\boxtimes_B Y$ to be the completion of $(X\otimes_\bbC Y) / N$
    under the norm
    $\|\xi\|_{X\boxtimes_B Y}^2:= \|\langle \xi|\xi\rangle_X\|_C$.
    
    The unit $A-A$ correspondence is given by ${}_AA_A$ with $\langle a_1|a_2\rangle_A :=a_1^*a_2$, and the unitors are given by the obvious unitary maps 
    $$
    \begin{aligned}
    \lambda^A_X:A\boxtimes_A X&\to X
    \\
    a\boxtimes \xi &\mapsto a\xi
    \end{aligned}
    \qquad\qquad
    \begin{aligned}
    \rho^B_X:X\boxtimes_B B &\to X
    \\
    \xi\boxtimes b &\mapsto \xi b
    \end{aligned}
    $$
    To construct the associator
    $\alpha_{X,Y,Z}:{}_A(X\boxtimes_B Y)\boxtimes_C Z_D \to {}_AX\boxtimes_B (Y\boxtimes_C Z)_D $, 
    we observe the map 
    $(\xi\otimes_\bbC \eta) \otimes_\bbC \zeta \mapsto \xi\otimes_\bbC (\eta \otimes_\bbC \zeta)$
    on the algebraic tensor product preserves $D$-valued inner products, and thus descends to a unitary isomorphism.
\end{ex}

We recall that generalizing the notion of Jones index for $\rm{II}_1$-factors \cite{MR696688}, the \textbf{right} (resp. \textbf{left}) \textbf{Watatani-index} for a Hilbert $A$-$B$-bimodule ${}_A X_B$ (i.e. $X$ has commuting  left-$A$ and a right-$B$ actions which are adjointable, and is simultaneously a left and right C*-correspondence with equivalent norms) \cite{MR996807, MR1624182}, provided it has a finite right (resp. left) {\bf Pimsner-Popa basis} $\{u_i\}_{i=1}^n\subset X$ ($\{v_j\}_{j=1}^m\subset X$): 
\begin{align}\label{eqn:WI}
    &\mathsf{Ind}_r(X):=\sum_{i=1}^n {}_A\langle u_i, u_i\rangle, && \mathsf{Ind}_\ell(X):=\sum_{j=1}^m \langle v_j\mid  v_j\rangle_A.
\end{align}
In case both left and right Watatani indices exist, we say that $X$ has \textbf{finite Watatani index} given by 
$$\mathsf{Ind}_W(X):=\mathsf{Ind}_r(X)\cdot \mathsf{Ind}_\ell(X),$$ 
which can be shown to be a positive element in the center of $A$. These indices are independent of the choices of Pimsner-Popa bases. 

\begin{defn}\label{defn:rigid}
    A C* 2-category $\cC$ is called \emph{rigid} if for every 1-morphism ${}_aX_b\in \cC(a\to b)$, there is a ${}_bX_a^\vee \in \cC(b\to a)$ together with maps $\ev{X} \in \cC({}_bX^\vee\boxtimes_a X_b \Rightarrow 1_b)$ and $\coev{X}\in \cC(1_a \Rightarrow {}_aX\boxtimes_b X^\vee_a)$ satisfying the zig-zag/snake equations, which are best depicted in the graphical calculus for 2-categories: 
    $$
    \tikzmath{
    \begin{scope}
    \clip[rounded corners=5pt] (-.7,-.6) rectangle (.7,.6);
    \fill[\AColor] (-.7,.6) -- (-.4,.6) -- (-.4,0) arc (-180:0:.2cm) arc (180:0:.2cm) -- (.4,-.6) -- (-.7,-.6);
    \fill[\BColor] (.7,.6) -- (-.4,.6) -- (-.4,0) arc (-180:0:.2cm) arc (180:0:.2cm) -- (.4,-.6) -- (.7,-.6);
    \end{scope}
    \draw (-.4,.6) -- (-.4,0) arc (-180:0:.2cm) arc (180:0:.2cm) -- (.4,-.6);
    }
    =
    \tikzmath{
    \begin{scope}
    \clip[rounded corners=5pt] (-.3,-.6) rectangle (.3,.6);
    \fill[\AColor] (-.3,-.6) rectangle (0,.6);
    \fill[\BColor] (.3,-.6) rectangle (0,.6);
    \end{scope}
    \draw (0,-.6) -- (0,.6);
    }
    \qquad\text{ and }\qquad
    \tikzmath{
    \begin{scope}
    \clip[rounded corners=5pt] (-.7,-.6) rectangle (.7,.6);
    \fill[\BColor] (-.7,-.6) -- (-.4,-.6) -- (-.4,0) arc (180:0:.2cm) arc (-180:0:.2cm) -- (.4,.6) -- (-.7,.6);
    \fill[\AColor] (.7,-.6) -- (-.4,-.6) -- (-.4,0) arc (180:0:.2cm) arc (-180:0:.2cm) -- (.4,.6) -- (.7,.6);
    \end{scope}
    \draw (-.4,-.6) -- (-.4,0) arc (180:0:.2cm) arc (-180:0:.2cm) -- (.4,.6);
    }
    =
    \tikzmath{
    \begin{scope}
    \clip[rounded corners=5pt] (-.3,-.6) rectangle (.3,.6);
    \fill[\BColor] (-.3,-.6) rectangle (0,.6);
    \fill[\AColor] (.3,-.6) rectangle (0,.6);
    \end{scope}
    \draw (0,-.6) -- (0,.6);
    }\,.
    $$
    Moreover, we assume for each 1-morphism ${}_aX_b\in \cC(a\to b)$, there is a \emph{predual} object ${}_b(X_\vee)_a \in \cC(b\to a)$ such that $(X_\vee)^\vee \cong X$ in $\cC(a\to b)$.
    Being rigid is a property of $\cC$, and not extra 
    structure.
    
    A \emph{unitary dual functor} on $\cC$ consists of a choice of dual $({}_bX^\vee_a, \ev{X}, \coev{X})$ for each 1-morphism ${}_aX_b\in \cC(a\to b)$ such that for each $f\in \cC({}_aX_b\Rightarrow {}_aY_b)$, we have $f^{\vee\dag}=f^{\dag\vee}$, where 
    $$
    f^\vee:=
    \tikzmath{
    \begin{scope}
    \clip[rounded corners=5pt] (-.9,-.8) rectangle (.9,.8);
    \fill[\AColor] (-.6,-.8) -- (-.6,.3) arc (180:0:.3cm) -- (0,-.3) arc (-180:0:.3cm) -- (.6,.8) -- (.9,.8) -- (.9,-.8);
    \fill[\BColor] (-.6,-.8) -- (-.6,.3) arc (180:0:.3cm) -- (0,-.3) arc (-180:0:.3cm) -- (.6,.8) -- (-.9,.8) -- (-.9,-.8);
    \end{scope}
    \draw (0,.3) arc (0:180:.3cm) -- (-.6,-.8);
    \draw (0,-.3) arc (-180:0:.3cm) -- (.6,.8);
    \roundNbox{fill=white}{(0,0)}{.3}{0}{0}{$f$}
    }.
    $$ 

    Whenever we have a unitary dual functor on $\cC,$ we define the involution: 
    $$f^*:= (f^\vee)^\dag = (f^\dag)^\vee.$$
\end{defn}

    In contrast to rigidity being a property, a unitary dual functor is extra structure. Notice we purposely suppressed the \emph{unitary tensorator} giving the coherence structure for the dual functor.
    We will also systematically suppress the \emph{canonical unitary pivotal structure} yielding a natural isomorphism $X\Rightarrow (X^\vee)^\vee$. These choices ensure that the {\bf left and right quantum dimensions} 
    $$d_L(X) = \ev{X}\circ\ev{X}^\dag\in \End(1_b),\  \text{ and }\  d_R(X) = \coev{X}^\dag\circ \coev{X}\in \End(1_a)$$ 
    are positive. 
    More generally, we have a {\bf left and right traces on}  $\End(1_b)$ and $\End(1_a)$ respectively, given by
    \begin{align*}
    \begin{aligned}
    \tr_L: \End(X)&\to \End(1_b)\\
    g &\mapsto 
    \tikzmath{
    \begin{scope}
    \clip[rounded corners=5pt] (-.9,-.8) rectangle (.7,.8);
    \fill[\BColor] (-.9,-.8) rectangle (.7,.8);
    \fill[\AColor] (0,.3) arc (0:180:.3cm) -- (-.6,-.3) arc(180:360:.3cm) -- (0,-.3);
    \end{scope}
    \draw (0,.3) arc (0:180:.3cm) -- (-.6,-.3) arc(180:360:.3cm) -- (0,-.3);
    \roundNbox{fill=white}{(0,0)}{.3}{0}{0}{$g$}
    }
    \end{aligned}
    \qquad \text{ and }\qquad
    \begin{aligned}
    \tr_R:\End(X)&\to \End(1_a)\\
        h&\mapsto
    \tikzmath{
    \begin{scope}
    \clip[rounded corners=5pt] (-.7,-.8) rectangle (.9,.8);
    \fill[\AColor] (-.7,-.8) rectangle (.9,.8);
    \fill[\BColor] (0,.3) arc (180:0:.3cm) -- (.6,-.3) arc(0:-180:.3cm) -- (0,-.3);
    \end{scope}
    \draw (0,.3) arc (180:0:.3cm) -- (.6,-.3) arc(0:-180:.3cm) -- (0,-.3);
    \roundNbox{fill=white}{(0,0)}{.3}{0}{0}{$h$}
    }
    \end{aligned}.
    \end{align*}
    In case $\tr_R=\tr_L$ for all objects, we say that {\bf the unitary dual functor is balanced}; as a consequence, the left and right quantum dimensions match.  
    This equality is the most meaningful whenever the endomorphism algebras $\End(1_a)$ are one-dimensional, and in such case we say the unitary dual structure {\bf is spherical}.  
    We reserve the notation $\overline{X}\cong X^\vee$ to refer to the choice of balanced dual for $X$, and similarly for the balanced dual of morphisms $\overline{f} = f^\vee$.    
    A more detailed discussion around unitary dual functors can be found in \cite{MR4133163}.

\begin{ex}\label{ex:ConjugateCorrespondence}
{\bf The conjugate correspondence.} Continuing the discussion after Example \ref{ex:rCorr}, given a right $B$-$A$ bimodule $Y$ of finite Watatani index $\mathsf{Ind}(Y)$ as in Equation \ref{eqn:WI}, we now describe the unitary dual $A$-$B$ bimodule ${_A}\overline{Y}_B$ together with its structure equations. In order to balance the resulting duality structure, we shall further assume that $Z(A)=\bbC1=Z(B),$ and denote $K:=\dfrac{\mathsf{Ind}_r(Y)}{\mathsf{Ind_\ell(Y)}}>0.$ 
Fix left and right Pimsner-Popa bases $\{v_j\}_1^m$ and $\{u_i\}_1^n$.
As a vector space, $\overline{Y}= \left\{\overline{\xi}\right\}_{\xi\in Y}$ is the conjugate space to $Y$, with the following $A$-$B$ bimodule structure: 
for all $a\in A,$ $b\in B$, and $\eta, \xi\in X$, 
\begin{align*}
    & a\rhd\overline{\xi}\lhd b:= \overline{b^*\rhd\xi\lhd a^*}, \text{ with inner products }\\
    \lip{A}{\overline{\eta}}{\overline{\xi}}:=&\rip{\eta}{\xi}{A}\quad \text{ and }\qquad   \rip{\overline{\eta}}{\overline{\xi}}{B}:= K^{-1/2}\lip{B}{\eta}{\xi}
\end{align*}
 
In coordinates, with respect to our chosen bases, and graphically, the associated duality maps for $Y$ are expressed as: 
    \begin{align}\label{eqn:DualityMaps}
        &\mathsf{ev}_Y:\overline{Y}\boxtimes_B Y\to A \qquad\qquad && \mathsf{coev}_Y:B\to Y\boxtimes_A\overline{Y}\\
        &\hspace{1.5cm}\overline{\eta}\boxtimes \xi\mapsto \langle\eta|\ \xi\rangle_A \qquad &&\hspace{1.5cm} b\mapsto  \left(\sum_{i=1}^nu_i\boxtimes\overline{u_i}\right)\lhd b\nonumber\\
        &\hspace{2.5cm}\tikzmath{
        \begin{scope}
        \clip[rounded corners=5pt] (-.7,-.8) rectangle (.7,.8);
        \fill[\AColor] (-.7,-.8) rectangle (.7,.8);
        \fill[\BColor] (.3,-.8)arc(0:180:.3cm) -- (0,-1);
    \end{scope}
            \draw(.3,-.8)arc(0:180:.3cm);    
            \draw[dotted](0,.-.4) --(0,.8);
        }
        &&\hspace{2.5cm}\tikzmath{
        \begin{scope}
        \clip[rounded corners=5pt] (-.7,-.8) rectangle (.7,.8);
        \fill[\BColor] (-.7,-.8) rectangle (.7,.8);
        \fill[\AColor] (.3,.8)arc(0:-180:.3cm) -- (0,1);
    \end{scope}
            \draw(.3,.8)arc(0:-180:.3cm);    
            \draw[dotted](0,.4) --(0,.-.8);
        }.\nonumber
    \end{align}
It is easy to see these satisfy the zig-zag equations from Definition \ref{defn:rigid}. 
A direct computation tells that
\begin{align*}
\begin{aligned}
    \ev{Y}^\dag (a)&= K^{1/2}\cdot a\sum_{j=1^m}\overline{v_j}\boxtimes v_j    
\end{aligned}
    \qquad \text{ and }\qquad
\begin{aligned}
    \coev{Y}^\dag(\eta\boxtimes\overline{\xi}) = \rip{\overline\eta}{\overline\xi}{B}, 
\end{aligned}
\end{align*}
which satisfy the balancing condition. 

We notice that if we remove the hypotheses on the centers of $A$ and $B$ and discard $K$, we still conclude $\overline{Y}$ is dual to $Y$, but this duality will not be necessarily balanced.     
\end{ex}

We define a \textbf{$\rm C^*$-tensor category} as a C* 2-category over a single object \cite[Remark 2.10]{MR4419534}. In practice, we forget about the object itself, but retain the data of its Hom 1-category and treat it as such. The following type of C*-tensor category plays a central role in our forthcoming considerations and examples: 

\noindent \textbf{Unitary tensor categories}. A \emph{unitary tensor category} (UTC) $\cC$ is a semisimple rigid $\rm C^*$-tensor category with simple unit object $1_{\cC}$. Rigidity means that every object has a unitary dual (see \cite[Definition 2.9]{MR4419534} and references therein), and simplicity of the tensor unit is the requirement that $\End_\cC(1_\cC)\cong\bbC$. If the isomorphism classes of simple objects in a unitary tensor category form a finite set, we say the category is a \emph{unitary fusion category}. 

\begin{ex}
    When $A=\bbC$, then $\rCorr(A\to A)$ simply becomes {\bf the UTC of Hilbert spaces} and bounded linear transformations $\Hilb.$ Moreover, the restriction to unitarily dualizable Hilbert spaces recovers the category of all finite dimensional Hilbert spaces $\fdHilb$, which is the familiar ambient where quantum information is customarily developed.     
\end{ex}

\begin{ex}[{\bf UTCs of finitely generated projective bimodules}]
    If $A$ is a unital $\rm C^*$-algebra, we define $\fgpBim(A)$ to be the full tensor subcategory of dualizable objects in $\rCorr(A\rightarrow A)$. 
\end{ex}

\begin{remark}
To justify the notation $\fgpBim(A)$, recall that by \cite{MR2085108}, any dualizable right $A$-$A$ correspondence admits a compatible structure of a left finitely generated projective Hilbert $A$-module. 
This is, if $X\in \rCorr(A\to A)$ is dualizable, using the evaluation and coevaluation maps, one can obtain a canonical left $A$-valued inner product, making $X$ into a left Hilbert $A$-module. It turns out that the left and right Hilbert $A$-module structure induce equivalent topologies on $X$ and the actions are bi-adjointable, and so $X$ is an $A$-$A$ Hilbert bimodule which is finitely generated and projective on the left and on the right.     
\end{remark}

\begin{defn}
    \textbf{Unitary tensor functors}. For unitary tensor categories $\cC$ and $\cD,$ a \emph{unitary tensor functor} is a triple $(F, F^1, F^2),$ where $F:\cC\to \cD$ is a $\dag$-functor, 
    $F^1:F(1_{\cC})\to 1_{\cD}$ is a chosen unitary isomorphism, and 
    $F^2 =\left\{F^2_{a,b}:F(a)\otimes F(b)\to F(a\otimes b)|\ a,b\in\cC \right\}$ is a natural unitary isomorphism satisfying the standard coherence axioms c.f. \cite[Chapter 2.4]{MR3242743}. Notice we don't require our unitary tensor functors be fully faithful. In this paper, we will assume our monoidal functors are \textit{strictly unital}, i.e. $F(1_{\cC})=1_{\cD}$ and $F^1=\id_{1_{\cD}}.$
\end{defn}

\subsection{Q-Systems}
\ \\
We shall now briefly recount the necessary background on Q-systems on C* 2-categories. Q-systems are categorified finite dimensional C*-algebras and therefore are appropriate categorical generalizations of the quantum sets considered in \cite{MR3849575} for the category of finite dimensional Hilbert spaces $\fdHilb$. For a detailed exposition on Q-systems we refer the reader to \cite[\S3]{MR4419534} and \cite{MR3308880}. Let $\cC$ denote an arbitrary C* 2-category.

\begin{defn}\label{defn:QSystem}
An algebra in $\cC$ is a 1-endomorphism $Q\in \cC(b\to b)$
with a multiplication $m\in\cC(Q\boxtimes_b Q\Rightarrow Q)$ and unit $i\in\cC(1_b\Rightarrow Q)$ 2-morphisms represented graphically by a trivalent and univalent vertices respectively, and $m^\dag$ and $i^\dag$ are given by the vertical reflections:
$$
\tikzmath{
\fill[\BColor, rounded corners=5pt] (-.3,0) rectangle (.9,.6);
\draw (0,0) arc (180:0:.3cm);
\draw (.3,.3) -- (.3,.6);
\filldraw (.3,.3) circle (.05cm);
}=m,
\qquad\qquad
\tikzmath{
\fill[\BColor, rounded corners=5pt] (-.3,0) rectangle (.9,-.6);
\draw (0,0) arc (-180:0:.3cm);
\draw (.3,-.3) -- (.3,-.6);
\filldraw (.3,-.3) circle (.05cm);
}=m^\dag,
\qquad\qquad
\tikzmath{
\fill[\BColor, rounded corners=5pt] (0,0) rectangle (.6,.6);
\draw (.3,.3) -- (.3,.6);
\filldraw (.3,.3) circle (.05cm);
}=i,
\qquad\qquad
\tikzmath{
\fill[\BColor, rounded corners=5pt] (0,0) rectangle (.6,-.6);
\draw (.3,-.3) -- (.3,-.6);
\filldraw (.3,-.3) circle (.05cm);
}=i^\dag.
$$
We call $Q$ a \emph{Q-system} if $(Q,m,i)$ satisfies:
\begin{enumerate}[label=(Q\arabic*)]
\item 
\label{Q:associativity}
(associativity)
$\tikzmath{
\fill[\BColor, rounded corners=5pt] (-.3,-.3) rectangle (1.2,.6);
\draw (0,-.3) -- (0,0) arc (180:0:.3cm);
\draw (.3,-.3) arc (180:0:.3cm);
\draw (.3,.3) -- (.3,.6);
\filldraw (.3,.3) circle (.05cm);
\filldraw (.6,0) circle (.05cm);
}
=
\tikzmath{
\fill[\BColor, rounded corners=5pt] (-.6,-.3) rectangle (.9,.6);
\draw (0,0) arc (180:0:.3cm) -- (.6,-.3);
\draw (-.3,-.3) arc (180:0:.3cm);
\draw (.3,.3) -- (.3,.6);
\filldraw (.3,.3) circle (.05cm);
\filldraw (0,0) circle (.05cm);
}$,
\item
\label{Q:unitality}
(unitality)
$\tikzmath{
\fill[\BColor, rounded corners=5pt] (-.3,-.3) rectangle (.9,.6);
\draw (0,-.1) -- (0,0) arc (180:0:.3cm) -- (.6,-.3);
\draw (.3,.3) -- (.3,.6);
\filldraw (.3,.3) circle (.05cm);
\filldraw (0,-.1) circle (.05cm);
}
=
\tikzmath{
\fill[\BColor, rounded corners=5pt ] (0,-.3) rectangle (.6,.6);
\draw (.3,-.3) -- (.3,.6);
}
=
\tikzmath{
\fill[\BColor, rounded corners=5pt] (-.3,-.3) rectangle (.9,.6);
\draw (0,-.3) -- (0,0) arc (180:0:.3cm) -- (.6,-.1);
\draw (.3,.3) -- (.3,.6);
\filldraw (.3,.3) circle (.05cm);
\filldraw (.6,-.1) circle (.05cm);
}$,
\item
\label{Q:Frobenius}
(Frobenius)
$
\tikzmath{
\fill[\BColor, rounded corners=5pt] (-.3,-.6) rectangle (1.5,.6);
\draw (0,-.6) -- (0,0) arc (180:0:.3cm) arc (-180:0:.3cm) -- (1.2,.6);
\draw (.3,.3) -- (.3,.6);
\draw (.9,-.3) -- (.9,-.6);
\filldraw (.3,.3) circle (.05cm);
\filldraw (.9,-.3) circle (.05cm);
}
=
\tikzmath{
\fill[\BColor, rounded corners=5pt] (-.3,0) rectangle (.9,1.2);
\draw (0,0) arc (180:0:.3cm);
\draw (0,1.2) arc (-180:0:.3cm);
\draw (.3,.3) -- (.3,.9);
\filldraw (.3,.3) circle (.05cm);
\filldraw (.3,.9) circle (.05cm);
}
=
\tikzmath{
\fill[\BColor, rounded corners=5pt] (-.3,.6) rectangle (1.5,-.6);
\draw (0,.6) -- (0,0) arc (-180:0:.3cm) arc (180:0:.3cm) -- (1.2,-.6);
\draw (.3,-.3) -- (.3,-.6);
\draw (.9,.3) -- (.9,.6);
\filldraw (.3,-.3) circle (.05cm);
\filldraw (.9,.3) circle (.05cm);
}
$,
\item
\label{Q:separable}
(separable)
$
\tikzmath{
\fill[\BColor, rounded corners=5pt] (-.3,0) rectangle (.9,1.2);
\draw (0,.6) arc (180:-180:.3cm);
\draw (.3,1.2) -- (.3,.9);
\draw (.3,0) -- (.3,.3);
\filldraw (.3,.3) circle (.05cm);
\filldraw (.3,.9) circle (.05cm);
}
=
\tikzmath{
\fill[\BColor, rounded corners=5pt ] (0,0) rectangle (.6,1.2);
\draw (.3,0) -- (.3,1.2);
}$.
\end{enumerate}
If $(Q,m,i)$ satisfies \ref{Q:associativity}, \ref{Q:unitality}, and \ref{Q:Frobenius}, we call $Q$ a \emph{$\rm C^*$ Frobenius algebra}.
\end{defn}

Notice that if $Q$ is a Q-system, then $Q$ is dualizable, self-dual in fact, with 
$$
\begin{aligned}
    \ev{Q}:=
\tikzmath{
\fill[\BColor, rounded corners=5pt] (-.5,-.1) rectangle (.5,.5);
\draw (-.2,.5) arc (-180:0:.2cm);
\draw (0,.3) -- (0,.1);
\filldraw (0,.3) circle (.05cm);
\filldraw (0,.1) circle (.05cm);
}
\end{aligned}
\qquad\text{ and }\qquad
\begin{aligned}
\coev{Q}
:=
\tikzmath{
\fill[\BColor, rounded corners=5pt] (-.5,-.5) rectangle (.5,.1);
\draw (-.2,-.5) arc (180:0:.2cm);
\draw (0,-.3) -- (0,-.1);
\filldraw (0,-.1) circle (.05cm);
\filldraw (0,-.3) circle (.05cm);
}    
\end{aligned}.
$$

Given a Q-system or $\rm C^*$-Frobenius algebra $(Q,m,i)$, we define
\begin{equation}\label{eq:dQ}
d_Q:=\tikzmath{
\fill[\BColor, rounded corners=5pt ] (0,0) rectangle (.6,1);
\draw (.3,.3) -- (.3,.7);
\filldraw (.3,.3) circle (.05cm);
\filldraw (.3,.7) circle (.05cm);
}
\in \End_\cC(1_b)^+.
\end{equation}
By \cite[Lem.~1.16]{MR2298822}, we have
$$
    \tikzmath{
    \fill[\BColor, rounded corners=5pt] (-.5,-.5) rectangle (.5,.5);
    \draw (-.2,-.5) arc (180:0:.2cm);
    \draw (-.2,.5) arc (-180:0:.2cm);
    \draw (0,-.3) -- (0,-.1);
    \draw (0,.3) -- (0,.1);
    \filldraw (0,.3) circle (.05cm);
    \filldraw (0,.1) circle (.05cm);
    \filldraw (0,-.1) circle (.05cm);
    \filldraw (0,-.3) circle (.05cm);
    }
    \leq 
    \tikzmath{
    \fill[\BColor, rounded corners=5pt] (-.9,-.5) rectangle (.5,.5);
    \draw (-.2,-.5) -- (-.2,.5);
    \draw (.2,-.5) -- (.2,.5);
    \draw (-.6,-.2) -- (-.6,.2);
    \filldraw (-.6,.2) circle (.05cm);
    \filldraw (-.6,-.2) circle (.05cm);
    }
    \leq
    \|d_Q\|\cdot
    \tikzmath{
    \fill[\BColor, rounded corners=5pt] (-.5,-.5) rectangle (.5,.5);
    \draw (-.2,-.5) -- (-.2,.5);
    \draw (.2,-.5) -- (.2,.5);
    }.
$$
Moreover, by \cite[Cor.~5.18]{MR2298822} we know that each C* Frobenius algebra in a tensor $\rm C^*$-category is equivalent to a Q-system satisfying \ref{Q:associativity}-\ref{Q:separable}.
We remark that $d_Q$ might fail to be invertible in this generality. 
In case $\End(1_b)$ is finite dimensional and all 1-morphisms in $\End(b)$ are dualizable, then there are various notions of {\bf quantum dimension} in $\End(b)$ \cite{MR1444286,MR3994584,MR4133163}. 
One may then consider \emph{normalized} Q-systems for which $d_Q$ is invertible and equals this quantum dimension.

\begin{ex}[{\bf Q-systems as pairs of pants}]\label{ex:Pants}
Let ${}_aX_b\in \cC(a\to b)$, where $\cC$ is a C* 2-category.
A {\bf unitarily separable left dual} for ${}_aX_b$ is a dual $({}_bX^{\vee}_a , \ev{X}, \coev{X})$ such that $\ev{X} \circ \ev{X}^{\dag} = \id_{1_b}$. (There is a similar notion of a unitarily separable right dual which we shall not use.) We will state whenever we are considering a unitarily separable left dual $X^\vee$, and keep reserving the notation $\overline{X}$ to denote the balanced dual. (Recall discussion following Definition \ref{defn:rigid}.)

Given a unitarily separable left dual for ${}_aX_b$, ${}_aX\boxtimes_b X^\vee_a \in \cC(a\to a)$ is a Q-system with multiplication and unit given by
$$
m:=
\tikzmath{
\begin{scope}
\clip[rounded corners=5pt] (-.7,0) rectangle (.7,.9);
\fill[\AColor] (-.7,0) rectangle (.7,.9);
\fill[\BColor] (-.4,0) -- (-.4,.2) .. controls ++(90:.2cm) and ++(270:.2cm) .. (-.1,.7) -- (-.1,.9) -- (.1,.9) -- (.1,.7)  .. controls ++(270:.2cm) and ++(90:.2cm) .. (.4,.2) -- (.4,0);
\fill[\AColor] (-.2,0) -- (-.2,.2) arc (180:0:.2cm) -- (.2,0);
\end{scope}
\draw (-.2,0) -- (-.2,.2) arc (180:0:.2cm) -- (.2,0);
\draw (-.4,0) -- (-.4,.2) .. controls ++(90:.2cm) and ++(270:.2cm) .. (-.1,.7) -- (-.1,.9);
\draw (.4,0) -- (.4,.2) .. controls ++(90:.2cm) and ++(270:.2cm) .. (.1,.7) -- (.1,.9);
}
=
\id_X\boxtimes \ev{X}\boxtimes \id_{X^\vee}
\qquad\text{ and }\qquad
i:=
\tikzmath{
\begin{scope}
\clip[rounded corners=5pt] (-.4,-.4) rectangle (.4,.5);
\fill[\AColor] (-.7,-.4) rectangle (.7,.5);
\fill[\BColor] (-.1,.5) -- (-.1,0) arc (-180:0:.1cm) -- (.1,.5);
\end{scope}
\draw (-.1,.5) -- (-.1,0) arc (-180:0:.1cm) -- (.1,.5);
}
=\coev{X}.
$$
\end{ex}

\begin{remark}
    It will often be the case in our applications that we can renormalize the unitarily separable left dual into the balanced unitary dual and vice versa. 
    Moreover, we will often explicitly describe our bimodules, their Pimsner-Popa bases and the duality equations we work with. 
    Thus, to simplify the notation, we will only refer to pairs of pants using conjugates (ie $\overline{X}X$) obfuscating the specific flavor of dual, and renormalizing our duality maps however we see fit. 
\end{remark}

\begin{ex}[{\bf Q-systems from finite groups}]\label{ex:Q-sysFiniteGroups}
    Let $G$ be a finite group and $\omega\in Z^{3}(G, \mathsf{U}(1))$ a unitary 3-cocycle on $G$, and consider the UTC of {\bf twisted $G$-graded Hilbert spaces} $\cC=\Hilb(G,\omega)$. Modulo isomorphism, the connected (i.e. $1_\cC$ appears exactly once in its semisimple decomposition) Q-systems in  are given by pairs $(H, \mu)$, where $H\le G$ and $\mu:H\times H\rightarrow \mathsf{U}(1)$ is a trivialization of $\omega|_{H}$. 
    Moreover, two $Q$-systems $(H, \mu)$ and $(K, \nu)$ are isomorphic if and only if $H=K$ and $\mu$ differs from $\nu$ by a coboundary. Moreover, $(H, \mu)$ and $(K, \nu)$ are Morita equivalent (i.e. isomorphic via a bimodule) if there exists a $g\in G$ such that $(H^{g},\mu^{g})$ is isomorphic to $(K, \nu),$ where $(\cdot)^{g}$ denotes the conjugation automorphism.
    See \cite{MR3933035}, \cite[Example 3.9]{MR4717816}. 
\end{ex}

\begin{ex}[{\bf Q-Systems from compact quantum groups}]\label{ex:Q-SysCompactQuantumGroups}
Q-systems also arise naturally when considering compact quantum groups. Let $\bbG$ be such, and by \cite[\S2.3]{MR3933035} and references therein, the finite dimensional $\bbG$-C*-algebras equipped with their canonical choice of faithful $\bbG$-invariant state correspond to Q-systems in $\Rep(\bbG)$, the unitary tensor category of finite dimensional representations of $\bbG$.
\end{ex}

Given a C* 2-category and a Q-system $Q\in\cC(b\to b),$ its endomorphisms $\cC(Q\Rightarrow Q)$, or simply $\End_\cC(Q)$, form a C*-algebra under the given $\dag$-structure from $\cC$ and its vertical composition $-\circ-$. However, using the structure maps from $Q,$ we can endow $\End_\cC(Q)$ with a second operation called the {\bf Schur product}. (Sometimes called in the literature Hadamard product or convolution product) Given $S,T\in \End_\cC(Q)$ we define their Schur product:
\begin{equation}\label{eqn:SchurProduct}
S\star T:= m\circ(S\boxtimes T)\circ m^\dag=\  
\tikzmath{
\begin{scope}
\clip[rounded corners=5pt] (-1.1,-1.5) rectangle (1.1,1.5);
\fill[\AColor] (-1.1,-1.5) rectangle (1.1,1.5);
\end{scope}
\draw (-.5,.3) arc (180:0:.5cm);
\draw (0,.8) -- (0,1.5);
\filldraw (0,.8) circle (.05cm);
\roundNbox{fill=white}{(-.5,0)}{.3}{0}{0}{$S$}
\roundNbox{fill=white}{(.5,0)}{.3}{0}{0}{$T$}
\draw (-.5,-.3) arc (180:360:.5cm);
\draw (0,-.8) -- (0,-1.5);
\filldraw (0,-.8) circle (.05cm);
}\ \in\End_\cC(Q).
\end{equation}

We say that $S\in\End_\cC(Q)$ is a {\bf Schur idempotent} if $S= S\star S.$ It is evident from the string diagrams that the Jones projection is a Schur idempotent:

\begin{align}\label{eqn:JonesShur}
e\star e\ =\ \tikzmath{
\begin{scope}
\clip[rounded corners=5pt] (-.7,-.9) rectangle (.7,.9);
\fill[\AColor] (-.7,-.9) rectangle (.7,.9);
    \fill[\BColor] (-.2,.2) arc (360:180:.1cm) -- (-.4,.2) .. controls ++(90:.2cm) and ++(270:.2cm) .. (-.1,.7) -- (-.1,.9) -- (.1,.9) --(.1,.7)  .. controls ++(270:.2cm) and ++(90:.2cm) .. (.4,.2) arc (360:180:.1cm);
\fill[\AColor] (-.2,0) -- (-.2,.2) arc (180:0:.2cm) -- (.2,0);
    \fill[\BColor] (-.2,-.2) arc (-360:-180:.1cm) -- (-.4,-.2) .. controls ++(-90:.2cm) and ++(-270:.2cm) .. (-.1,-.7) -- (-.1,-.9) -- (.1,-.9) --(.1,-.7)  .. controls ++(-270:.2cm) and ++(-90:.2cm) .. (.4,-.2) arc (-360:-180:.1cm);
\fill[\AColor] (-.2,0) -- (-.2,-.2) arc (-180:-0:.2cm) -- (.2,0);    
\end{scope}
\draw  (-.2,.2) arc (180:0:.2cm) -- (.2,.2);
\draw  (-.4,.2) .. controls ++(90:.2cm) and ++(270:.2cm) .. (-.1,.7) -- (-.1,.9);
\draw (-.4,.2) arc (180:360: .1cm);
\draw (.4,.2) arc (360:180: .1cm);
\draw (.4,.2) .. controls ++(90:.2cm) and ++(270:.2cm) .. (.1,.7) -- (.1,.9);
\draw  (-.2,-.2) arc (-180:-0:.2cm) -- (.2,-.2);
\draw  (-.4,-.2) .. controls ++(-90:.2cm) and ++(-270:.2cm) .. (-.1,-.7) -- (-.1,-.9);
\draw (-.4,-.2) arc (-180:-360: .1cm);
\draw (.4,-.2) arc (-360:-180: .1cm);
\draw (.4,-.2) .. controls ++(-90:.2cm) and ++(-270:.2cm) .. (.1,-.7) -- (.1,-.9);
}
\ =
\tikzmath{
\begin{scope}
\clip[rounded corners=5pt] (-.5,-.9) rectangle (.5,.9);
\fill[\AColor] (-.7,-1) rectangle (.7,1);
\fill[\BColor] (-.2,.9) -- (-.2,.7) arc (-180:0:.2cm) -- (.2,.9);
\fill[\BColor] (-.2,-.9) -- (-.2,-.7) arc (-180:-360:.2cm) -- (.2,-.9);
\end{scope}
\draw (-.2,.9) -- (-.2,.7) arc (-180:0:.2cm) -- (.2,.9);
\draw (-.2,-.9) -- (-.2,-.7) arc (-180:-360:.2cm) -- (.2,-.9);
}
\ = e.
\end{align}

\begin{ex}\label{ex:ClassicalSet}
    If $\cC=\fdHilb$ viewed as a rigid C* 2-category over a single object, for $N\in\bbN$, we consider the commutative Q-system $Q=\bbC^N$ with $\{\delta_\ell\}_1^N,$  the canonical basis, and usual multiplication $m(\delta_\ell\otimes \delta_{k})= \delta_{\ell=k}\delta_\ell$ and unit $1\mapsto N^{-1}\sum_1^N \delta_\ell$, then one can check directly that the comultiplication is given by the linear extension of $m^\dag(\delta_\ell) = \delta_\ell\otimes \delta_\ell,$ and the counit  by $i^\dag(\delta_\ell)=N^{-1}1_{\bbC}$ for all $\ell,k=1,2,...,N.$ 
    Moreover, if $S,T\in\End(\bbC^N)\cong M_N$, then the Schur product $S\star T\in M_N$ becomes entry-wise multiplication 
    $$(S\star T)_{\ell k}= S_{\ell k}\cdot T_{\ell k}\qquad\forall \ell,k=1,2,...,N.$$
\end{ex}

\subsection{Subfactors and the standard invariant}\label{sec:SubfacIndexStdInv}
\ \\
In this section we will summarize the constructions and tools from subfactors that we will make use of. A subfactor is a unital inclusion of von Neumann algebras with trivial center, and a popular useful reference is the textbook \cite{MR1473221}. 

We consider a unital inclusion of C*-algebras $A\subset B$ equipped with a faithful conditional expectation (i.e. an idempotent ucp $A$-$A$ bimodular map) $E:B\twoheadrightarrow A.$ Let $\cB$ be the completion of $B$ as a right Hilbert $A$-module with inner product: $$\langle b_1\mid  b_2\rangle_A := E(b_1^* b_2).$$ 
There is a canonical faithful embedding of $B$ into $\cB,$ 
whose image we denote by $B\Omega\subset\cB.$ The norm induced by the $A$-valued inner product on $B\Omega$ is denoted on elements by  $\|b\Omega\|_A:=\|\langle b\Omega\mid b\Omega\rangle_A\|^{1/2}.$ 
We notice that $\cB\in \rCorr(B\to A)$ as $B$ acts on $\cB$ on the left by adjointable operators by (extending) the multiplication map from $B\Omega$ to all $\cB$. 

A sufficient condition so that $\cB\cong B\Omega$ as right C* $A$-modules, ie $\|\cdot\|_B\cong\|\cdot\|_A$ on $B$, so it is $\|\cdot\|_A$-complete, is that there is some constant $c\geq 1$ such that
$\|\cdot\|_B\leq c\|\cdot\|_A,$ since $\|\cdot\|_A\leq \|\cdot\|_B$ is automatic by the positivity of $E$. 
We say that $E$ has finite \textbf{Pimsner--Popa index}, if such $c$ giving the norm comparison exists. 
In this case, we define the value of the Pimsner--Popa index as 
$$\textsf{Ind}_p(E)= \textsf{inf}\{c\geq 1|\ (cE-\id)\geq 0\}.$$
If no such $c\geq 1$ exists, we say the Pimsner--Popa index for the inclusion is infinite.  \cite[Definition 2]{MR1642530}, \cite{MR860811}.

In \cite[Proposition 2.1.5]{MR996807}, Watatani establishes that finiteness of the Pimsner--Popa index together with unitality of the right $A$-compact operators $\cK(B_A)$ (i.e. the \textbf{reduced C*-basic construction} for $A\subset B$ satisfies $\cK(B_A) = \End(B_A)$) are equivalent to the finiteness of the Watatani index. In this work we shall mostly be interested in subfactors and C*-inclusions of finite Jones-Watatani index, so we will refer to the Hilbert C*-module $\cB$ simply as $B$. 

\begin{defn}\label{defn:JonesProjBasicConst}
    Given a unital inclusion of C*-algebras $A\overset{E}{\subset}B$ of finite Watatani index $[B:A]:=\mathsf{Ind}_W(E)$, the {\bf Jones projection} $e\in\End({}_AB_A)$ is the orthogonal projection onto $A$ given by $e(b\Omega) = E(b)\Omega.$ Here, for $b\in B$ we write $b\Omega$ to distinguish elements in $B$ viewed as a bimodule from operators $b\in B.$ The existence of $e$ is guaranteed since $A_A$ is an fgp (=finitely generated projective) submodule of $B_A$ and it is thus orthogonally complemented (cf. \cite[Lemma 1.9]{2023arXiv230505072H}).  
    Then $e$ implements $E$ on $B$ (as a multiplication operator on $B\Omega$) as: 
    $$ E(b)\circ e= e\circ b \circ e.$$
    Up to a normalization factor, the Jones projection corresponds to the diagram 
    $$[B:A]\cdot e=\tikzmath{
        \begin{scope}
            \clip[rounded corners=5pt] (-.5,0) rectangle (.5,1);
            \fill[\AColor] (-.5,0) rectangle (.5,1);
            \fill[\BColor] (.3,0)arc(0:180:.3) -- (.3,0);
            \fill[\BColor] (-.3,1)arc(180:360:.3) -- (-.3,1);
        \end{scope}
        \draw(.3,0)arc(0:180:.3);
        \draw[dotted](0,.3)--(0,.7);
        \draw(-.3,1) arc(180:360:.3);
    }\in\End({}_AB_A).$$
    In this context, we have that ${\rm{C}}^*(B, e)\subseteq \End(B_A)$ is the reduced C*-basic construction (or {\bf the Jones basic construction} in the jargon of subfactors), and so ${\rm{C}^*}(B, e)= \End(B_A) = \cK(B_A)$ \cite{MR996807}. 
\end{defn}
We remark that these considerations allow for many flexible ways to present the reduced C*-basic construction:
\begin{equation}\label{eqn:BasicConst}
    \End(B_A)\cong \mathsf{span}\left\{beb'|\ b,b'\in B\right\} \cong \mathsf{span}\left\{|b\rangle\langle b'|\ |\  b,b'\in B\right\}\cong\cK(B_A).
\end{equation}

We shall now demonstrate how to present {\bf the standard invariant of an inclusion} $(A\overset{E}{\subset}B)$ (i.e. the lattice of higher relative commutants) \cite{MR1334479} in the language of C* 2-categories and Q-systems \cite{MR1966524}.
\begin{construction}\label{construction:HigherRelComs}
    The basic ingredient is {\bf the Jones tower}. Given a unital inclusion of C*-algebras $A\subset B$ with a faithful conditional expectation $E:B\twoheadrightarrow A$ onto $A$ with finite Watatani index $\mathsf{Ind}_W(E)\in Z(A)$, we perform the reduced C*-basic construction to obtain a unital C*-algebra $B_1:=\End(B_A).$ The left action of $B_0:=B$ on $B_A$ by multiplication then yields an inclusion $B_0\subset B_1$, which is moreover equipped with a faithful expectation  $E_1:B_1\twoheadrightarrow B_0$ given by the linear extension of $beb'\mapsto bb'$, using Equation \eqref{eqn:BasicConst}. 
    Crucially, $B\overset{E_1}{\subset} B_1$ is a unital inclusion, where $E_1$, the \emph{dual conditional expectation}, has the same Watatani index as the initial $E$ \cite[Proposition 2.3.4]{MR996807}.
    Moreover, we can identify 
    $$
    B\boxtimes _A B\cong B_1
    $$ 
    as C*-algebras via the isomorphism $B\otimes_A B\ni x\boxtimes y \mapsto (x\rhd)\circ e_A\circ(y\rhd)\in \End(B_A).$
    Here, $x\rhd$ is the operator on $B$ by left multiplication by $x$ and $e_A$ is the Jones projection down to $A$ \cite[Lemma 2.2.2]{MR996807}.
    
    Therefore, the reduced C*-basic construction of a finite index inclusion can be iterated, obtaining a chain of inclusions, all with fixed index $\rm{Ind}_W(E).$ Relabeling $B_{-1}:= A$ and $E_0:=E$ we obtain the Jones tower: 
    $$A=B_{-1}\overset{E_0}{\subset} B_0 \overset{E_1}{\subset} B_1 \overset{E_2}{\subset} B_2 \subset \cdots \subset B_{n-1}\overset{E_n}{\subset} B_n \subset \cdots.$$
    By the remarks in the previous paragraph, for $n\in \bbN$, we can identify 
    $$
    B_n\cong (B\boxtimes_A B)^{\boxtimes_B n}
    $$
    as C*-algebras. 

    We are now ready to introduce the standard invariant of $A\overset{E}{\subset}B$ as {\bf the lattice of higher relative commutants}: 
    \begin{equation}\label{eq:HigherRelComms}
    \begin{tikzcd}[column sep=0em, row sep=-.03em]
     A'\cap A\ 
     \subset\ 
     &  A'\cap B\  
     \subset\ 
     & A'\cap B_1 \ 
     \subset\
     & \cdots\ 
     \subset\ 
     & A'\cap B_n\ \subset\ \cdots\\
     {}  
     &\hspace{-.6cm}\cup  
     &\hspace{-.7cm}\cup 
     &{}
     &\hspace{-1.1cm}\cup\\
     &  B'\cap B\ 
     \subset\ 
     &  B'\cap B_1\ 
     \subset\    
     & \cdots\
     \subset
     & B'\cap B_n \subset\hdots
     \end{tikzcd}.
    \end{equation}
    Which is a lattice of inclusions of finite dimensional C*-algebras, provided that $A'\cap A$ is finite dimensional. When $A'\cap A= \bbC1$, as one moves right in either row, the dimension cannot increase faster than a factor of $[B:A]=\mathsf{Ind}_W(E)$ \cite[Lemma 4.6.2 ii]{MR999799}.

    We shall now recast this lattice in the language of C* 2-categories and Q-systems.  From our starting finite-index inclusion $A\subset B$, we obtain a dualizable $B$-$A$ correspondence ${}_BB_A$ with dual ${}_A\overline{B}_B$ given by the conjugate correspondence from Example \ref{ex:ConjugateCorrespondence}, which explicitly gives its dual. 
    We now consider $\cC_{A\subset B},$ which is the C* 2-category obtained by taking subobjects, direct sums, duals and Connes fusion tensor product. In $\cC:=\cC_{A\subset B}$ We find the pair of pants (recall Example \ref{ex:Pants}) Q-system given by $Q:= {}_A\overline{B}\boxtimes B_A$, which is manifestly split by $X:={}_BB_A$ and $\overline{X}:= {}_A\overline{B}_B.$

    The lattice of higher relative commutants \eqref{eq:HigherRelComms} then becomes
     \begin{equation}\label{eq:HigherRelCommsEnd}
    \begin{tikzcd}[column sep=0em, row sep=-.03em]
     \End_\cC({}_AA_A)
     \subset 
     &  \End_\cC({}_A\overline{X}_A)  
     \subset 
     & \End_\cC({}_A\overline{X}\boxtimes_BX_A)  
     \subset
     & \End_\cC({}_A\overline{X}\boxtimes_BX\boxtimes_A\overline{X}_B) 
     \subset \cdots\\
     {}  
     &\hspace{-.6cm}\cup  
     &\hspace{-.7cm}\cup 
     &\hspace{-1.1cm}\cup\\
     &  \End_\cC({}_BB_B)\    
     \subset 
     &\hspace{.5cm}  \End_\cC({}_BX_A) \ \hspace{.5cm} 
     \subset    
     &\hspace{.5cm} \End_\cC({}_BX\boxtimes_A\overline{X}_B)\ \hspace{.5cm}
     \subset \cdots
     \end{tikzcd}.
    \end{equation}
Here, the horizontal inclusions are given by composing with right tensoring by $X$ and $\overline{X}$ in an alternating pattern, and the vertical inclusions are given by composing with left tensoring by $\overline{X}.$ 

Finally, since $X$ is dualizable, we can rely on {\bf Frobenius reciprocity} realize these endomorphism algebras as alternating products of $X$ and its conjugate. For instance, by repeatedly \emph{bending strings up}, one gets:
$$\End_{\cC}(\overline{X} X) = \cC(\overline{X} X\to \overline{X} X)\cong \cC(\overline{X} \to \overline{X} X\overline{X}) \cong \cC(A\to \overline{X} X\overline{X}X)\cong \overline{X} X\overline{X} X.$$
\end{construction}
In Section \ref{sec:QFourierTransf} we will provide explicit computations of the higher relative commutants in various specific examples.

\subsection{The quantum Fourier transform}\label{sec:QFourierTransf}
Every unitary tensor category $\cC$ embeds in the UTC of fgp bimodules over a $\rm{II}_1$-factor \cite{BHP12} and by \cite{MR4139893}, $\cC$ also embeds into $\fgpBim(A)$ for some unital separable C*algebra $A$ with a unique trace. 
That is, we can always represent UTCs as fgp bimodules over some operator algebra. 
We therefore give the following definition 
\begin{defn}\label{defn:ConcreteCat}
    We say that a C* 2-category $\cC$ is {\bf concrete} if there exists a $\dag$ 2-functor $$F:\cC\to \rCorr.$$ We call $F$ a {\bf generalized fiber functor}. Since there are often many fiber functors out of a given $\cC$, when we refer to a concrete C* 2-category, we often denote it by $(\cC, F)$ whenever necessary. 
\end{defn}
This is a generalization of the \emph{concrete tensor categories}, which is a tensor category admitting a \emph{fiber functor} (ie monoidal) into $\fdHilb$ considered by \cite{MR3849575}. 
If $\cC$ is a UTC, the unitary tensor functors $\cC\to\fgpBim(A)$ turn $\cC$ into a concrete C* 2-category, by considering $\cC$ as a 2-category over a single object in the usual way. 
    
The main advantage of having a concrete C* 2-category is that $\rCorr$ (Example \ref{ex:rCorr}) is \emph{Q-system complete} \cite{MR4419534}. 
This can be quickly interpreted by saying that every Q-system in $\rCorr$ unitarily splits; i.e. is equivalent to a pair of pants (recall Example \ref{ex:Pants}). We record this in the following remark: 

\begin{remark}\label{rmk:ModelQSys}
    If $Q\in\rCorr(A\to A)$ is a concrete Q-system over the unital C*-algebra $A$, then (\cite[Example 3.12, Corollary 4.17]{MR4419534}) there exists $B\in\rCorr$ and a dualizable $B$-$A$ correspondence ${}_BX_A(={_B}B_A)$ with unitary dual ${}_A\overline{X}_{B}(={_A}B_B)$ such that $Q\cong \overline{X}\boxtimes_{B} X$ as $Q$-systems. 
    
    In short, we can always concretely represent $Q$ as a shaded double strand $\overline{X}\boxtimes X$ with the pair of pants (co-)multiplication. 
    Moreover, we can always afford to model any Q-system in any concrete C* 2-category in terms of a unital inclusion of C*-algebras $A\subset B$ equipped with a faithful surjective cp $A-A$ bimodular map $E: B\to A$ such that $B_A$ is finitely generated projective.
\end{remark}

The pair of pants perspective of Q-systems affords the use of the quantum Fourier transform (cf \cite[\S3]{MR1442497} \cite{Ocn90, MR4236188}), which will become an indispensable tool in the sequel. 
\begin{defn}\label{defn:FT}
    Let $A,B$ be unital C*-algebras, and let $X\in\rCorr(B\to A)$ be  a dualizable correspondence. 
    The, {\bf quantum Fourier transform} is the linear isomorphism  corresponding to a graphical rotation: 
    \begin{align*}
        \cF_2: \rCorr_{B-B}(X\boxtimes_A\overline{X}) &\to\rCorr_{A-A}(\overline{X}\boxtimes_B X)\\
        \tikzmath{
            \begin{scope}
                \clip[rounded corners=5pt] (-.6,1) rectangle (.6,-1);
                \fill[\BColor](-.6,-1) rectangle (.6,1);
                \fill[\AColor] (-.2,1) -- (-.2,-1) -- (.2,-1) -- (.2,1);
            \end{scope}
            \draw (-.2,.3) -- (-.2,1);
            \draw (.2,.3) --(.2,1);
            \roundNbox{fill=white}{(0,0)}{.3}{.1}{.1}{$T$}
            \draw (-.2,-.3) -- (-.2,-1);
            \draw(.2,-.3) --(.2,-1);
        }&\mapsto\quad
        \tikzmath{
        \begin{scope}
            \clip[rounded corners=5pt] (-1.1,-1) rectangle (1.1,1);
            \fill[\AColor] (-1.1,-1) rectangle (1.1,1);
            \fill[\BColor] (-.8,1)--(-.8,-.3) arc(180:360:.3) --(-.2,1);
            \fill[\BColor] (.2,-1) -- (.2,.3) arc(180:0:.3) -- (.8,-1);
        \end{scope}
            \draw(-.2,-.3) arc(360:180:.3) -- (-.8,1);
            \draw(.2,-.3) -- (.2,-1);
            \roundNbox{fill=white}{(0,0)}{.3}{.1}{.1}{$T$}
            \draw(.2,.3) arc(180:0:.3) -- (.8,-1);
            \draw (-.2,.3) -- (-.2,1);
        }=: \hat{T}.
    \end{align*}
\end{defn}
In other words, $$\hat{T} = (\id_{\overline{X}X}\boxtimes \ev{X})\circ(\id_{\overline{X}}\boxtimes T\boxtimes \id_X)\circ (\ev{X}^\dag\boxtimes\id_{\overline{X}X}).$$
Observe that $\cF_2$ maps the identity to the Jones projection, and vice versa as a direct consequence of the zig-zag equations. 
We will later interpret this fact after we have established certain properties of the Fourier transform. We denote $\cF^{-1}(S)$ by $\check{S},$ and will often drop the $2$ from the subindex when the context is clear.  

\begin{lem}\label{lem:FTIso}
    The Fourier transform 
    $$\cF_2: \left(\rCorr_{B-B}(X\boxtimes_A\overline{X}),\ \star \right) \to\left(\rCorr_{A-A}(\overline{X}\boxtimes_B X),\ \circ \right)$$ is a $\bbC$-linear isomorphism that transforms the Schur product $\star$ to the categorical vertical composition $\circ$. 
    More precisely, 
    $$\cF_2(T_1\star T_2) = \cF_2(T_1)\circ\cF_2(T_2).$$
    Additionally, $\cF_2$ induces a bijection between the positive cones of $(\rCorr_{B-B}(\overline{X}\boxtimes X), \star)$ and $(\rCorr_{A-A}(X\boxtimes \overline{X}, \circ).$
    
    Moreover, if we swap the operations of composition and Schur product, we obtain that 
    $$\cF_2: \left(\rCorr_{B-B}(X\boxtimes_A\overline{X}),\ \circ,\ \dag \right) \to\left(\rCorr_{A-A}(\overline{X}\boxtimes_B X),\ \star,\ * \right)$$
    is a $*$-algebra anti-isomorphism, satisfying 
    $$\cF_2(T^\dag)= \cF_2(T)^*.$$ 
    Furthermore, $\cF_2$ remains an isomorphism when restricted to $\circ$-idempotents. 
    Recall that $T^\dag$ is the adjoint map to $T$ in $\rCorr$, and if $S$ is a morphism, then $S^*= \overline{S}^\dag = \overline{(S)^\dag}$ is the dual of $S$ followed by the $\dag$ in $\cC$ (c.f. Definition \ref{defn:rigid} and the discussion following it).
\end{lem}
\begin{proof}
    From a direct graphical argument, it is clear that $\cF_2$ is an invertible $\bbC$-linear map between finite dimensional vector spaces and that it is exchanges the two operations as in the statements above. 
    
    The restriction of $\cF_2$ to $\star$-idempotents remains one-to-one so it suffices to show it is onto $\circ$-idempotents. 
    Let $S=S\circ S\in \rCorr(\overline{X}\boxtimes X).$ Then, by the first statement its inverse transform $\check S$ is a $\star$-idempotent. Since $\cF$ is an isomorphism, then $\cF_2(\check{S})= S$ implies that the restriction of $\cF_2$ to $\star$-idempotents is also surjective onto $\circ$-idempotents. 
    That $\cF_2$ is an isomorphism from $\circ$-idempotents onto $\star$-idempotents is proven similarly. 
    
    The identity $\cF_2(T^*)= \cF_2(T)^\dag$ is readily verified by a direct graphical computation. 
    
    Finally, the bijection between $\star$-positives and $\circ$-positives is established similarly. 
\end{proof}

\subsection{The Fourier transform and the dual of a finite group}

We will concretely interpret the results from Section \ref{sec:QFourierTransf} through the  motivating example of a finite group, while also making direct contact with aspects of the theory of subfactors outlined in Section \ref{sec:SubfacIndexStdInv}. 

\begin{ex}[{\bf The Fourier transform and the dual of a finite group}]\label{ex:GroupHyperfiniteCommutants} (cf \cite[Example from p. 39]{Ocn90}) 
    Let $\Gamma$ be a finite group. 
    Then $\Gamma$ has an essentially unique outer action on $\cR,$ the hyperfinite $\rm{II}_1$-factor \cite{MR587749}.  
    We can then form a subfactor by taking a crossed product by this action, and build the Jones tower (Construction \ref{construction:HigherRelComs}):
    $$\overbrace{\cR}^{N}\ \overset{E}{\subset}\ \overbrace{\cR\rtimes\Gamma}^{M}\  \overset{E_1}{\subset}\ \overbrace{\underbrace{M_\Gamma(\cR)}_{\cong \cR\otimes\cB(\ell^2\Gamma)}}^{M_1}\ \overset{E_2}{\subset}\  \overbrace{\underbrace{M_\Gamma(\cR\rtimes \Gamma)}_{\cong (\cR\rtimes\Gamma)\otimes\cB(\ell^2\Gamma)}}^{M_2}\  \subseteq ...,$$
    where each step is obtained as a crossed-product by (co)-actions of $\Gamma$, but also matches the Jones basic construction. 
    Notice the first Jones projection $e\in M_\Gamma(\cR)$ is given by $(e)_{ij} = \delta_{i=j=\id},$ so that the index of the inclusions is given by $\tr(e)^{-1}=|\Gamma|<\infty.$
    The beginning of the lattice of higher relative commutants is given by :
    \begin{equation}\label{ex:GpStdInv}
    \begin{tikzcd}[column sep=0em]
     \overbrace{\bbC1}^{N'\cap N}
     \hookrightarrow
     &  \overbrace{\bbC1}^{N'\cap M}
     \hookrightarrow
     & \overbrace{\ell^\infty (\Gamma)}^{N'\cap {M_1}} 
     \hookrightarrow
     & \overbrace{\mathsf{Mat}_\Gamma(\bbC)}^{N'\cap {M_2}}  
     \hookrightarrow
     & \hdots\\
     {}
     &  \underbrace{\bbC1}_{M'\cap M}
     \arrow[hook]{u}
     \hookrightarrow
     &  \underbrace{\bbC1}_{M'\cap {M_1}}
     \arrow[hook]{u}
     \hookrightarrow
     \arrow["\cF_1"']{ul}
     & \underbrace{\bbC[\Gamma]}_{M'\cap {M_2}} 
     \arrow[hook]{u}
     \arrow["\cF_2"']{ul}
     \hookrightarrow
     & \hdots
     \end{tikzcd}.
    \end{equation}
    That $N'\cap M\cong\ell^\infty(\Gamma)$ follows from the fact that $$\cR\rtimes\Gamma = \bigoplus_{g\in\Gamma}\ {}_g\cR$$ as a $\Gamma$-graded C*-algebra, where for each $g\in\Gamma,$ ${}_g\cR$ is the $\cR$-$\cR$ bimodule with diagonal action on the right and left action $r\rhd_g \xi= (g\cdot r)\xi$ for $r,\xi\in \cR.$ 
    We let $\{u_g\}_{g\in\Gamma}\subset M$ be a collection of unitaries implementing the action; this is, 
    $$\forall x\in N,\quad \forall g\in \Gamma,\qquad g\cdot x = u_g x u_{g}^*.$$
    Along these lines, outerness boils down to
    $$\forall g\in \Gamma, \forall x\in N\qquad  xu_g=u_gx\implies g=\id_\Gamma.$$
    Then, $\End_{\cR-\cR}(\cR\rtimes\Gamma)$ is a direct sum of finite dimensional C*-algebras indexed over the group elements, and so we write ${_g}\cR=u_g\cR$. 
    Furthermore, the size of each matrix block equals the multiplicity of each $g$-twisted bimodule ${_g}\cR$ inside the crossed-product; namely: one. 
    Concretely, if $P_g$ is the projection on $\cR\rtimes_\alpha\Gamma$ onto ${_g}\cR$ for $g\in\Gamma,$ the relative commutant $N'\cap M$ is thus 
    $$\bigoplus_{g\in\Gamma} P_g\bbC = \End_{\cR-\cR}(\cR\rtimes_\alpha\Gamma)\cong\ell^\infty(\Gamma).$$ 
    Alternatively, one could use outerness of the action to directly compute the relative commutant $\cR'\cap \End(\cR\rtimes\Gamma_\cR),$ where $\cR$ acts on the left by multiplication by elements in ${}_e\cR$.

    We now use the Fourier transform to compute the relative commutant $M'\cap M_2.$ Indeed, Lemma \ref{lem:FTIso} tells us that $(\End({_M}{M_1}_{M}), \circ)\cong (\End({_N}{M}_N), \star),$ and since the latter is spanned by the projections $\{P_g\}_g,$ it suffices to check that $P_a\star P_b = P_{ab}$ for all $a,b\in \Gamma$ to conclude that $M'\cap M_2\cong \bbC[\Gamma^{\op}]$ with the $-\circ-$ composition.
    Indeed, the Schur product is implemented as 
    $$P_a\star P_b = m\circ (P_a\boxtimes_N P_b)\circ m^\dag,$$
    where $m^\dag: M\to M\boxtimes_N M$ is the comultiplication map which is given by (cf \cite[Lemma 3.11]{MR4419534})
    $$m^\dag(x)\mapsto \sum_{g\in \Gamma} u_g\boxtimes_Au_g^*\cdot x.$$
    Finally, a direct computation will show that 
    $$P_a\star P_b(u_c)= \delta_{ab=c}u_{ab}= P_{ab}(u_c) \qquad \forall a,b,c\in \Gamma,$$
    and so we conclude that $M'\cap M_2\cong \bbC[\Gamma^{\op}]$ with product by $-\circ-$ composition. 
    
    Therefore, {\bf all Schur idempotents arising from the subfactor $N\subset M$ live in $N'\cap M_1\cong\ell^\infty(\Gamma)$, and correspond with the $\circ$-idempotents in $\bbC[\Gamma]$ via the Fourier transform}. 
    We will reconnect with this setting in Example \ref{ex:GammaGraphs}.     
\end{ex}

\section{Categorified sets and graphs}\label{sec:IntroCatQuSetGraph}

We shall introduce the notions of sets and graphs internal to a C* 2-category or, perhaps more suggestively, lay out a {\bf $\cC$-equivariant theory of sets and graphs}. 
The idea is to provide a unifying framework to encompass the various formulations of classical and quantum sets and graphs, as well as their  morphisms, naturally appearing the study of non-commutative objects and in quantum information theory.
A substantial advantage of our proposed framework is its enormous  flexibility, as it affords passage from the abstract to the concrete by incorporating tools from C*-algebras and subfactors, providing interesting reinterpretations of the latter, and new constructions for quantum sets and graphs.

Quantum sets and graphs have been studied and formalized by several authors in different contexts \cite{2009arXiv0906.2527D, MR2908249, MR3549479, MR3849575, MR3907958, MR3926289, MR4302212,MR4481115, MR4514486, MR4555986, MR4706978} and known to be equivalent under reasonable constraints. 
In its most elementary form, a finite quantum set is a finite dimensional C*-algebra $(B,\ m:B\otimes B\to B)$ equipped with a \textbf{delta form} $\psi:B\to \bbC$, which is a faithful state satisfying $m \circ m^\dag = \delta^2 \id_B.$ Here, $m^\dag$ is the adjoint of the linear map $m:L^2(B,\psi)\otimes L^2(B,\psi)\to L^2(B,\psi)$ in the familiar GNS representation with respect to $\psi$. 
Musto, Reutter and Verdon \cite{MR3849575} encoded this idea via \emph{special symmetric dagger Frobenius algebras} in $\fdHilb$, and following Gelfand duality obtained a $\dag$ 2-category $\QSet$ of finite quantum sets, quantum functions and intertwiners.

Naturally, a finite quantum graph is defined as a pair $\cG=(B, \hat T),$ where $B$ is the \emph{finite quantum set of vertices}, and $\hat T$ is its \emph{quantum adjacency matrix}. 
The key observation here is that this formulation is entirely internal to $\fdHilb$, the smallest unitary tensor category consisting of finite-dimensional Hilbert spaces. And, that all this can as well be modeled in a subfactor or a C* 2-category without referring to Hilbert spaces, therefore incorporating new quantum sets and graphs into the same theory. 

We notice that this categorification mindset liberates us from finite-dimensional constraints while keeping a tight algebraic grip on the operators we will use. The ongoing theme here will be to replace classical notions such as dimension and orthonormal bases by quantum analogues such as the Jones-Watatani index/quantum dimension and Pimsner-Popa bases.  
This framework can be quite nicely described in terms of unitary tensor categories and Q-systems, which possess familiar diagrammatic calculi, and such that computations are still carried out in finite dimensional settings.

\begin{defn}\label{defn:CatQuSet}
A categorified finite quantum set in a C* 2-category $\cC$ is a Q-system $(Q, m, i)$ in $\cC.$ 
\end{defn}

If $\cC = \fdHilb$ (considered as a 2-category over one object), then the above definition in the symmetric case recovers the usual notion of finite quantum set modelled on finite dimensional C*-algebras  (equipped with a $\delta$-form). In this context, classical sets correspond to finite dimensional commutative C*-algebras; i.e. $\bbC^N$ together with the state coming from the uniform measure.

We shall now introduce categorified graphs. In analogy with the classical case, a simple finite directed graph $\cG$ with vertex set $V$ is specified by a $V\times V$ matrix $\hat T$ with entries in $\{0,1\}.$ 
We can capture this latter condition using the Schur product of matrices of the same size: $(S\star S')_{u,v} = S_{u,v}S'_{u,v}$; ie the entry-wise product. The adjacency matrix $\hat T$ is therefore a \textbf{Schur idempotent}. That is, $\hat T\star \hat T= \hat T.$ 
\begin{defn}\label{defn:QuGraph}
    Given a C* 2-category $\cC$, and a categorified quantum set $Q\in\cC(a\to a)$ for some object $a\in \cC$. A {\bf $\cC$-equivariant graph} $\cG$ is a pair $(Q, \hat{T}),$ where $\hat{T}\in \End_{\cC}(Q)$ is (a positive invertible scalar multiple of) a Schur idempotent. 
    Here, scalar multiple is with regards to  $Z(\End_\cC(Q)),$ which is an abelian C*-algebra.  
    We say that $\cG=(Q,\hat T)$ is
    \begin{itemize}
        \item {\bf Self-adjoint} if $\hat{T}= \hat{T}^\dag$.
        \item {\bf Real} (ie $*$-preserving) if $\hat{T}^*=\hat{T}$. 
        \item {\bf Undirected}  if $\hat{T}$ is both real and self-adjoint. 
        \item {\bf Reflexive} if $\hat{T}\star\id_Q = z\id_Q$ for some $z\in Z(\End_{\cC}(Q))$, and irreflexive if  $\hat{T}\star\id_Q = 0.$
    \end{itemize}
    Recall that the dagger $\hat{T}^\dag$ corresponds to the dagger involution in $\cC,$ while the $*$-operation is given by:
    $$T^* = (T^\dag)^{\wedge} = (T^\wedge)^\dag,$$
    where $T^\wedge=\overline{T}$ is the unitary dual of $T$ from Definition \ref{defn:rigid}.

    As has been explained by the aforementioned authors, $\hat{T}$ being real corresponds to $\hat T$ being $*$-preserving, and reflexivity corresponds to a graph having self-loops at every vertex in the classical case. 
    Notice that in this generality, since we are not assuming $\cC$ comes with a braiding (and such choice of structure is not canonical at times), we do not define an undirected graph via $\hat{T}\circ \mathsf{braid} = \hat{T},$ but instead chose the definition above, since it already extends the classical examples. It would be interesting to contrasts these notions in braided tensor categories. 
    The passage between a reflexive and an irreflexive graph is seamless and can be done by adding or substracting $\id_Q$ to the adjacency data. 

We will drop $\cC$ from the notation whenever specifying the category is not important, and simply make reference to categorified graphs, denoted $\cG=(Q,\hat T).$ 
\end{defn}

Particular examples of categorified graphs have appeared in the subfactors literature under the name of \emph{biprojections}. 
These are bimodule endomorphisms of a Q-system that are idempotent vertically and, up to a scalar factor, also horizontally. 
We recall their formal definition:
\begin{defn}\label{defn:biprojection}
    Let $\cC$ be a C* 2-category and $Q\in\cC(b\to b)$ be a Q-system of the form $Q\cong X\boxtimes_a \overline{X}$. 
    An element $p\in \End_{\cC}(Q)$ is called a {\bf biprojection} if and only if $p=p\circ p$ is a $\circ$-idempotent and $\cF_2(p)=z\cdot \cF_2(p)\star\cF_2(p)$ for some $z\in Z(\End_\cC(Q))_+^\times.$ (Recall Definition \ref{defn:FT}.)
\end{defn}

Biprojections are known to correspond to intermediate subfactors, and so we record this useful fact in the following example:
\begin{ex}[{\bf Biprojections from intermediate subfactors}]\label{ex:BiprojSubalgs} (\cite{MR1262294, MR1950890})
    Consider a unital inclusion of C*-algebras $A\subset B$ with a conditional expectation $E:B\twoheadrightarrow A$ of finite Jones/Watatani index. 
    From the graphical computation in Equation \ref{eqn:JonesShur}, we know that the Jones projection from Definition \ref{defn:JonesProjBasicConst} satisfies  $e\star e=e,$ so $e$ yields an explicit (canonical!) example of a Schur idempotent, and thus a biprojection. 
    Similarly, so is $\id_B$, since ${}_AB_A$ is a Q-system. 

    In fact, more is true: there is a bijective correspondences between biprojections $p\in\End_{A-A}(B)$ and unital intermediate C*-algebras $D$  (i.e. $A\subset D\subset B)$ with common unit. \cite[Theorem, p. 37]{Ocn90}
\end{ex}
As we shall see in Example \ref{ex:BiProjectionsCombinatorially}, biprojections carry graph-theoretic content, which classically encodes the adjacency matrix of disjoint unions of cliques. 

\medskip

Consistent with \cite[\S2.3]{MR4514486}, we say {\bf the set of quantum vertices of} $\cG=(Q, \hat T)$ consists of the irreducible 1-morphisms appearing in the semisimple decomposition of $Q$ as a 1-morphism in $\cC$ (recorded with multiplicities), whenever this is defined. (This is always the case whenever $\cC$ is a C* 2-category with finite dimensional endomorphism C*-algebras.) We express this as $Q\cong \oplus_i Q_i,$ where the $Q_i$ are irreducibles. 
In this context, {\bf the number of vertices} of $\cG$ is given by $$\#V_\cG:=d_Q=i^\dag\circ i= \sum_i d_{Q_i},$$
the quantum dimension of $Q.$ (See Definition \ref{defn:QSystem}.) Furthermore, it was already noted by Gromada that the number of irreducibles $Q_i$ appearing in the semisimple decomposition of $Q$ is not invariant under quantum isomorphisms even when $\cC=\fdHilb.$ 
We now highlight a fact about real Schur idempotents that is already established for the usual $\fdHilb$-equivariant graphs.
\begin{prop}
    Let $\cC\subset \rCorr$ be a concrete C* 2-category and let $\cG=(Q\cong {}_AB_A, \hat{T}\in \End_{A-A}(B))$ be a $\cC$-equivariant concrete graph. Then $\hat{T}$ is real if and only if $\hat T$ is positive (as an operator on C*-algebras i.e. $\hat{T}(b^*b)\geq 0\ \forall b\in B$) if and only if $\hat{T}$ is completely positive.
\end{prop}
\begin{proof}
    The proof is identical to \cite[Proposition 2.23]{MR4481115} but instead using the graphical calculus of $\cC$ and a Pimsner-Popa basis for $Q$, and so we omit it.
\end{proof}

We shall now discuss the space of edges for a general $\cC$-equivariant graph. 
A direct diagrammatic computation shows that the morphism:
\begin{align}\label{eqn:EdgeProjector}
    E_\cG:=\tikzmath{
    \begin{scope}
        \clip[rounded corners=5pt] (-1.1,-1.2) rectangle (1.1,1.2);
        \fill[\AColor] (-1.1,-1.2) rectangle (1.1,1.2);
    \end{scope}
    \draw (-.3,.6)--(-.3,1.2) node[above]{$\scriptstyle Q$};
    \draw (0,.3) arc (0:180:.3cm) -- (-.6,-1.2) node[below]{$\scriptstyle Q$};
    \filldraw (-.3,.6) circle (.05cm);
    \filldraw (.3,-.6) circle (.05cm);
    \draw (.3,-.6)--(.3,-1.2) node[below]{$\scriptstyle Q$};
    \draw (0,-.3) arc (180:360:.3cm) -- (.6,1.2) node[above]{$\scriptstyle Q$};
    \roundNbox{fill=white}{(0,0)}{.3}{.1}{0}{$\hat{T}$}  
    }
    =(m_Q\otimes \id_Q)\circ (\id_Q\otimes \hat{T}\otimes \id_Q) \circ (\id_Q\otimes m_Q^\dag)\in \cC(Q\otimes Q).
\end{align}
defines a $\circ$-idempotent, called {\bf the edge projector}, whose range corresponds to the {\bf edge space} of $\cG.$ (cf \cite[Definition 15]{MR3549479} and  \cite{MR3849575, MR4555986}) 
We now define {\bf the quantum edges of $\cG$} as the irreducible objects appearing in the semi-simple decomposition of the underlying 1-morphism determined by $E_\cG$. 
We can furthermore {\bf count the number of edges} as 
$$\#(E_\cG):= (\tr_L)\circ (\id_Q\otimes \tr_R) (E_\cG)= \iota_Q^\dag \circ \hat{T}\circ \iota_Q.$$
In the case of classical graphs (Example \ref{ex:ClassicalGraphs}), these constructions readily recover the usual meaning of an edge and count them appropriately. 
In a followup work \cite{BHPII}, we will revisit this construction in the context of concrete C* 2-categories (i.e. $\rCorr$) and equipped with the property of \emph{Q-system splittings} and the quantum Fourier Transform. But for now, this will suffice for the purposes of this section. 

In order to convince the reader that quantum $\cC$-equivariant graphs can be described quite concretely, we dedicate Section \ref{sec:Examples} to write down explicit examples, both known and new, of categorified graphs.
In a future work, we will  return to developing their abstract theory in detail, showing $\cC$-equivariant graphs assemble into a 2-category, which will open the door to developing an equivariant classification and connections into quantum information theory and invariants for operator algebras. It is timely, however, at this stage to introduce {\bf classical homomorphisms between quantum $\cC$-equivariant graphs}:
\begin{defn}\label{defn:ClassicalHomcCGraph}
    Let $\cC$ be a C* 2-category, and let $\cG=(Q, \hat{T})$ and $\cG'=(Q', \hat{T'})$ be $\cC$-equivariant graphs. A classical homomorphism between $\cG$ and $\cG'$ is a unital $*$-cohomomorphism of the underlying Q-systems $f\in\cC(Q\Rightarrow Q')$ preserving the adjacency relation; this is, 
    \begin{align*}
    (f\otimes f)\circ E_\cG=
    \tikzmath{
    \begin{scope}
        \clip[rounded corners=5pt] (-.6,-.8) rectangle (1.1,1.6);
        \fill[\AColor] (-.6,-.8) rectangle (1.1, 1.6);
    \end{scope}
        \roundNbox{fill=white}{(-.2,.9)}{.25}{0}{0}{$f$}
        \roundNbox{fill=white}{(.7,.9)}{.25}{0}{0}{$f$}
        \draw(-.2,1.15)--(-.2,1.6) node[above]{$\scriptstyle Q'$};
        \draw(.7,1.15)--(.7,1.6) node[above]{$\scriptstyle Q'$};
        \draw (-.2,.3)--(-.2,.65);
        \draw (.7,.3)--(.7,.65);
        \roundNbox{fill=white}{(0,0)}{.3}{0}{.5}{$E_\cG$}
        \draw (-.2,-.3)--(-.2,-.8) node[below]{$\scriptstyle Q$};
        \draw (.7,-.3)--(.7,-.8) node[below]{$\scriptstyle Q$};
    } = 
    \tikzmath{
    \begin{scope}
        \clip[rounded corners=5pt] (-.6,-.8) rectangle (1.1,2.6);
        \fill[\AColor] (-.6,-.8) rectangle (1.1, 2.6);
    \end{scope}    
        \roundNbox{fill=white}{(0,1.8)}{.3}{0}{.5}{$E_{\cG'}$}
        \draw (-.2,2.1)--(-.2,2.6);
        \draw (.7,2.1)--(.7,2.6);
        \roundNbox{fill=white}{(-.2,.9)}{.25}{0}{0}{$f$}
        \roundNbox{fill=white}{(.7,.9)}{.25}{0}{0}{$f$}
        \draw(-.2,1.15)--(-.2,1.5);
        \draw(.7,1.15)--(.7,1.5);
        \draw (-.2,.3)--(-.2,.65);
        \draw (.7,.3)--(.7,.65);
        \roundNbox{fill=white}{(0,0)}{.3}{0}{.5}{$E_\cG$}
        \draw (-.2,-.3)--(-.2,-.8);
        \draw (.7,-.3)--(.7,-.8);
    }    
    = E_{\cG'}\circ (f\otimes f)\circ E_\cG.
    \end{align*}
\end{defn}
The conditions on the function $f$ above come from the theory of quantum sets and quantum functions spelled by Musto, Reutter and Verdon. 
They characterized bijections of a quantum set  as unital $*$-co-homomorphisms of the underlying algebras which are also co-unital homomorphisms. 
In this case, we say $f$ is a {\bf classical graph isomorphism of categorified graphs} if $f:Q\to Q'$ is a bijection of the underlying quantum sets, and preserves the adjacency data: 
$$
f\circ \hat{T} = \hat{T'}\circ f.
$$ 
Details are spelled-out in \cite[\S 5]{MR3849575} and we omit them here.

\section{All categorified graphs}\label{sec:Examples}
In this section we show how to explicitly construct all finite classical and quantum graphs using the tools described in the previous sections. 
Moreover, we will describe the categorified graphs supported on the $A_n$-subfactors.

We first explain how to canonically endow any categorical quantum set $Q$ in a C* 2-category $\cC$ with the adjacency of a complete graph. Every other finite graph with vertex quantum set $Q$ will consequentially be a sub-graph in a precise sense specified later. 

Every Q-system in $\rCorr$ (respectively $\vNAlg$, the W* 2-category of von Neumann algebras and W*-correspondences \cite[\S2]{MR4419534}) arises from a finite Watatani-index unital inclusion $A\overset{E}{\subset} B$ of C*-algebras (resp. $\rm{II}_1$-factors), where $E:B\to A$ is a finite-index conditional expectation (faithful ucp map). \cite[\S3, Corollary 4.17]{MR4419534}. 
Recalling Definition \ref{defn:ConcreteCat} and the discussion preceding it, we can always model $Q$ as ${_A}B_A\cong {}_AB\boxtimes_BB_A\in\fgpBim(A),$ the concrete Q-system arising from some finite index unital inclusion of C*-algebras, turning $B$ into its underlying quantum set with coefficients in $A$. 
In this section, we shall adopt the perspective of finite-index inclusions $A\overset{E}{\subset} B$ as initial data.

\begin{construction}[{\bf Complete and trivial categorified graphs over quantum sets}]\label{const:CompleteTrivQuGraph}
    Let $A\overset{E}{\subset} B$ be a finite-index unital inclusion of C*-algebras, and let $X:={}_BB_A$ be a right $B$-$A$ correspondence, with $\langle b_1|\ b_2 \rangle_A =  E(b_1^*b_2).$ (Recall the ingredients from Construction \ref{construction:HigherRelComs} and the discussion at the beginning of Section \ref{sec:SubfacIndexStdInv}.)
    Let's assume $E$ has finite Watatani index $[B:A]\in Z(A)_+^\times$, then $X$ is dualizable with ${}_A\overline{X}_B = {}_AB_B.$ 
    Notice that $\overline{X}\boxtimes_B X = {}_A(B\boxtimes_B B)_A \cong {}_AB_A\cong Q$ recovers the underlying categorical quantum set, while ${}_B(X\boxtimes_A\overline{X})_B= {}_B(B\boxtimes_A B)_B\cong B_1$ yields the reduced C*-basic construction (Definition \ref{defn:JonesProjBasicConst}). 
    
    We now endow $Q$ with the structure of the \textbf{complete categorified graph} canonically. 
    Recall that the Jones projection (Definiton \ref{defn:JonesProjBasicConst}) is a multiple of a Schur idempotent (Example \ref{ex:BiprojSubalgs}), and so according to Definition \ref{defn:QuGraph}, 
    $$\cK_Q:= (Q={}_AB_A,\ e)$$ 
    is a categorified graph with vertex quantum set $B$ with scalars $A$, whose corresponding {\bf quantum edges} are $B\boxtimes_A \overline{B},$ which exactly matches the Jones basic construction! 

    Similarly, endowing $Q$ with the Schur idempotent $\id_{B}$ yields {\bf the trivial categorified graph} 
    $$\cT_Q:=\left( Q,\ \id_{B}\right).$$
    Notice that, by analogy, $C^*(B, \id_B) = B$ yields simply $B$ as the space of quantum edges.
\end{construction}

We will fully justify our nomenclature and analogies from Construction \ref{const:CompleteTrivQuGraph} by Example \ref{ex:ClassicalGraphs}, but first we reconnect with the quantum graphs which are currently receiving attention in the literature. 

\begin{ex}[{\bf The usual quantum graphs}]\label{ex:UsualQuGraphs}
    The quantum graphs considered by Musto-Reutter-Verdon are defined over quantum spaces supported over finite dimensional C*-algebras \cite{MR3849575}. In our language, these corresponds to Q-systems in $\fdHilb,$ arising from inclusions of the form $\bbC=A\overset{\psi}{\subset} B,$ where $B$ is a finite dimensional C*-algebra and $\psi$ a faithful state. Within this framework, any finite quantum graph defined over $\fdHilb$ will be given by a subspace of the Jones basic construction $B_1\cong B\otimes B$ in the spirit of \cite[\S 7]{MR3849575}. 
\end{ex}

We now explicitly describe categorified graphs which are not necessarily complete or trivial. 
Recall Lemma \ref{lem:FTIso}, whose relevance to us lies in the fact that the linear isomorphism $\cF_2$ allows us to {\bf obtain all possible Schur idempotents on a given categorical quantum set.}  
All is needed is to compute two steps of the towers of relative commutants $\{A'\cap B_i\}_{i\geq-1}$ and $\{B'\cap B_i\}_{i\geq 0},$ where $B_{-1}:=A$ and $A\subset B\subset B_1\subset B_2\subset\cdots$ is the Jones Tower.

With the notation and conventions from Construction \ref{const:CompleteTrivQuGraph}, with $X:={}_AB_B$ and $\overline X:={}_BB_A$ then $\cF_2$ becomes a linear isomorphism  
$$\cF_2:B'\cap B_2\to A'\cap B_1.$$ 
Recall from Construction \ref{construction:HigherRelComs} that we can identify $B_1\cong B\boxtimes_A B \cong C^*\langle B, e\rangle\cong \End(B_A)\subset \End(B)$. 
So that $B_2\cong \End((B_{1})_B)$.  
Thus, decomposing the finite dimensional C*-algebra $B'\cap B_2\cong \End({}_B(B_{1})_B)$ as a multi-matrix algebra will yield all its $\circ$-idempotents (including orthogonal projections) by direct inspection (modulo Murray-von Neumann equivalence). These, in turn  correspond to all the Schur idempotents inside $\End_{B-B}(B\boxtimes_A B)$ via $\cF_2.$

We will now do this explicitly for finite classical graphs.
\begin{ex}[{\bf All finite directed simple classical graphs}]\label{ex:ClassicalGraphs} 
    We model a classical set with $N$ indistinguishable points with the commutative C*-algebra $\bbC^N$ equipped with the tracial state extending uniform probability distribution. 
    Recalling the ingredients from Construction \ref{construction:HigherRelComs}, this is readily encoded by the inclusion 
    $$\bbC I_N=A\overset{\tr}{\subset}B= \bbC^N.$$ 
    Here, the Jones projection $e\in \End({_\bbC}\bbC^N_\bbC)\cong M_N(\bbC)$ implementing the trace $E=\tr$ has the following matricial expression with respect to the basis $\{N^{1/2}\delta_\ell\}_{\ell=1}^N$ in $\bbC^N:$
    $$Ne=\begin{pmatrix}
        1 & 1 & \hdots & 1\\
        1 & 1 & \hdots & 1\\
        \vdots & \vdots& \ddots & \vdots\\
        1 & 1 & \hdots & 1
    \end{pmatrix}.$$
    More generally, but still classically, we can think of the inclusion 
    $$\bbC I_N=A\overset{\omega}{\subset} B= \bbC^N,$$
    where $\omega$ is a faithful state on $B$ determined by the weights $\{q_\ell\}_1^N\subset(0,1)$ such that $\sum_1^Nq_\ell =1.$ In this way we regard an $N$ element classical set whose points have different sizes or charges. Whenever it is convenient, we think of $B$ as the diagonal subalgebra of $M_N(\bbC).$
        
    Let $X$ be the $B$-$A$ bimodule from Construction \ref{construction:HigherRelComs}, and $\overline{X}$ its dual. It is useful to recall the Q-system structure of $\bbC^N$ from Example \ref{ex:ClassicalSet}. This bimodule $X$ is equipped with the right $A$-linear inner product $\rip{\eta}{\xi}{A} = \omega(\eta^*\xi).$ 
    We consider the right Pimsner-Popa basis 
    $$\left\{u_\ell = q_\ell^{-1/2}\delta_\ell\right\}_{\ell=1}^N,$$
    which is readily seen to be orthonormal with respect to this product. Moreover, the index of the inclusion $A\overset{\omega}{\subset}B$ is given by 
    $$[B:A]= \sum_1^N u_\ell u_\ell^*= \mathsf{diag}\{q_1^{-1}, q_2^{-1}, \hdots, q_N^{-1}\}\in (A'\cap B)^+.$$ 
    We notice that in the tracial case, when $q_\ell =1/N$ for every $\ell=1,...,N$ the index becomes the scalar matrix $NI_N,$ which recovers the dimension of $B.$ 
    we can express the right $B$-valued inner product on $X$ as $\lip{B}{\eta}{\xi}=[B:A]^{-1/2}\eta\xi^*,$ and with it we succinctly write the necessary duality equations for $X$ to compute the quantum Fourier transform:
    \begin{align*}
    \begin{aligned}
        \ev{X}:\overline{X}\boxtimes X&\to \bbC\\
        \overline{\eta}\boxtimes\xi&\mapsto \rip{\eta}{\xi}{A}=\omega(\eta^*\xi)
    \end{aligned}
    \qquad\text{ and }\qquad
    \begin{aligned}
        \ev{X}^\dag:\bbC&\to \overline{X}\boxtimes X\\
        z&\mapsto [B:A]^{1/2}\sum_{\ell=1}^Nq_\ell \overline{u_\ell}\boxtimes u_\ell
    \end{aligned}.
    \end{align*}
    Here, to ease the notation, we write $\otimes$ to denote the relative tensor product $\boxtimes_\bbC,$ and $\boxtimes$ instead of  $\boxtimes_{\bbC^N}$ and we shall continue to do so throughout this example. 
    (The reason the index appears in the left $B$-valued inner product is so that the resulting duality equations give the unitary balanced dual. A computation reveals that $\ev{X}\circ\ev{X}^\dag=[B:A]^{1/2}$.)

    We shall now describe the Jones basic construction and the double tower of higher relative commutants up to level $2$: 
    \begin{equation}\label{ex:ClassicalStdInv}
    \begin{tikzcd}[column sep=0em]
         &  \overbrace{\cong \bbC^N}^{\overbrace{ \End\left({_\bbC}\bbC^N{_{\bbC^N}}\right) }^{\cong A'\cap B }}
         \hookrightarrow
         & \overbrace{\cong M_N(\bbC)\boxtimes I_N}^{\overbrace{\End\left({_\bbC}\bbC^N\boxtimes \bbC^N{_\bbC}\right)}^{\cong A'\cap {B_1}}} 
         \hookrightarrow
         & \overbrace{\cong   M_N(\bbC)\otimes \bbC^N}^{\overbrace{\End\left( {_\bbC}\bbC^N\boxtimes \bbC^N\otimes \bbC^N{_{\bbC^N}}\right)}^{\cong A'\cap {B_2}}}\\
    %     {}
         &  \underbrace{\cong\bbC^N}_{\underbrace{\End\left({_{\bbC^N}}\bbC^N{_{\bbC^N}}\right)}_{\cong B'\cap B}}
         \arrow[hook]{u}
         \hookrightarrow
         &  \underbrace{\cong \bbC^N\boxtimes  I_N}_{\underbrace{\End\left({_{\bbC^N}}\bbC^N\boxtimes \bbC^N{_{\bbC^N}}\right)}_{\cong B'\cap {B_1}}}
         \arrow[hook]{u}
         \hookrightarrow
         & \underbrace{\cong \bbC^N\boxtimes I_N\otimes\bbC^N }_{\underbrace{\End\left({_{\bbC^N}}\bbC^N\boxtimes \bbC^N\otimes \bbC^N{_{\bbC^N}} \right)}_{\cong B'\cap {B_2}} }
         \arrow[hook]{u}
         \arrow["\cF_2"']{ul}
    \end{tikzcd}.
    \end{equation}
    Here, we have highlighted the quantum Fourier transform $\cF_2$, mapping idempotents $T\circ T=T\in \End\left({_{\bbC^N}}\bbC^N\boxtimes_{\bbC^N}\bbC^N\boxtimes_{\bbC}\bbC^N_{\bbC^N} \right)$  into Schur idempotents $\hat T\star \hat T =\hat T\in \End\left({_\bbC}\bbC^N\boxtimes_{\bbC^N}\bbC^N_\bbC\right).$

    We shall now compute the quantum Fourier transform of an arbitrary $T\in \End{X\otimes\overline{X}}.$ 
    Writing $T$ with respect to our Pimsner-Popa basis gives 
    $$T(\overline{u_\gamma}\otimes u_\zeta)= \sum_{\lambda,\mu=1}^{N} T_{\gamma\zeta}^{\lambda\mu} u_{\lambda}\otimes\overline{u_{\mu}}$$
    At once, Definition \ref{defn:FT} yields
    \begin{equation}\label{eqn:QFTCN}
        \hat{T}(\overline{u_\alpha}\boxtimes u_\eta) = [B:A]^{1/2}\sum_{\ell, \lambda=1}^N q_\ell T^{\lambda\mu}_{\ell\alpha}\overline{u_\ell}\boxtimes u_\lambda
    \end{equation}
    On an elementary compact operator $T= |u_{a'}\rangle\langle u_{c'}|\otimes |\overline{u_{a}}\rangle\langle\overline{u_{c}}|$ we get: 
    \begin{align}\label{eqn:QFTCNCompact}
        \cF_2(|u_{a'}\rangle\langle u_{c'}|\otimes |\overline{u_{a}}\rangle\langle\overline{u_{c}}|)(\overline{u_\alpha}\boxtimes u_\eta) &= [B:A]^{1/2}\sum_{\ell, \lambda=1}^N q_\ell (|u_{a'}\rangle\langle u_{c'}|\otimes |\overline{u_{a}}\rangle\langle\overline{u_{c}}|)_{\alpha\eta}^{\beta\kappa} \overline{u_\ell}\boxtimes u_\lambda\\
        &=[B:A]^{1/2}\sum_{\ell, \lambda=1}^N q_\ell \delta_{a'=\lambda}\delta_{c'=\ell}\delta_{a=\eta}\delta_{c=\alpha} \overline{u_\ell}\boxtimes u_\lambda\nonumber\\
        &= [B:A]^{1/2}q_{c'}\delta_{a=\eta}\delta_{c=\alpha} \overline{u_{c'}}\boxtimes u_{a'}\nonumber
    \end{align}
    
    Now, by Lemma \ref{lem:FTIso} we know that $\cF_2$ restricts to a bijection between $\circ$-idempotents in $B'\cap B_2$ onto $\star$-idempotents in $A'\cap B_1.$ 
    This therefore recovers all finite simple directed classical graphs, since $\star$-idempotents in $A'\cap B_1 \cong M_N(\bbC)$ are precisely all $\{0,1\}$ matrices.  That is, all possible adjacency matrices, thus exhausting all finite classical simple directed graphs. Notice that while all our $\star$-idempotents here are necessarily real in this case, and so correspond to self-adjoint $\circ$-idempotents. 
    
    We shall now explicitly obtain all these $N\times N$ $\{0,1\}$-matrices  from Equation \eqref{eqn:QFTCNCompact}. First, we consider the unitary isomorphism $\lambda_X:\overline{X}\boxtimes X\to B$ given by $\overline{\eta}\boxtimes \xi \mapsto \eta^*\xi,$ whose adjoint is given by $\xi\mapsto \overline{I_N}\boxtimes\xi,$ so our adjacency matrices are indexed over $\{1,2,\hdots, N\}$. 
    From Equation \ref{eqn:QFTCNCompact} we have that, 
    $$
    \left[\lambda_X\circ\cF_2\big(|u_{a'}\rangle\langle u_{c'}|\otimes |\overline{u_{a}}\rangle\langle\overline{u_{c}}|\big)\circ\lambda_X^\dag\right]_{\eta b} = [B:A]^{1/2}q_{c'}^{1/2}q_{c}^{1/2}\delta_{a=\eta}\delta_{c'=a}\delta_{b=a'}. 
    $$
    Applied to the $\circ$-idempotent $|u_{a'}\rangle\langle u_{a'}|\otimes |\overline{u_{c}}\rangle\langle\overline{u_{c}}|$ composed with $\mathsf{ad}(\lambda_X)$ we obtain:
    \begin{align}
        \left[\lambda_X\circ \cF_2\big( |u_{a'}\rangle\langle u_{a'}|\otimes |\overline{u_{c}}\rangle\langle\overline{u_{c}}|\big)\circ\lambda_X^\dag\right]_{\eta b} &= [B:A]^{1/2}q_{a'}^{1/2}q_c^{1/2}\delta_{c=\eta}\delta_{c'=a'}\delta_{b=a'}\\
        & = \sqrt{\dfrac{q_c}{q_{a'}}}E_{a'c},\nonumber
    \end{align}
    where $E_{ij}$ is the $N\times N$ matrix whose only nonzero entry occurs at its $ij$-th entry. 
    This is on the nose a scalar multiple of an edge indicator. 

    We close this example exhibiting the Jones projection in the setting of a generic (non-uniform measure) on $N$ points with respect to our Pimsner-Popa basis:    $$\left[\lambda_X\circ\Big(\cF_2([B:A]^{-1/2}\ev{X}^\dag\circ\ev{X}) \Big)\circ\lambda_X^\dag \right]=\begin{pmatrix}
        q_1 & \sqrt{q_2q_1} & \cdots & \sqrt{q_Nq_1}\\
        \sqrt{q_1q_2} & q_2 & \cdots & \sqrt{q_Nq_2}\\
        \vdots & \vdots & {} & \vdots\\
        \sqrt{q_{1}q_{N-1}} & \sqrt{q_{2}q_{N-1}} & \cdots & \sqrt{q_Nq_{N-1}}\\
        \sqrt{q_{1}q_{N}} & \sqrt{q_{2}q_{N}} & \cdots & q_N\\
    \end{pmatrix}.$$
   Note that we recover the \emph{all ones} matrix from the Jones projection we saw above, in case all the $q_\ell=1/N.$ 
\end{ex}
\medskip 

Example \ref{ex:ClassicalGraphs} fully justifies the nomenclature and ideas introduced in Construction \ref{const:CompleteTrivQuGraph} regarding the complete and trivial categorified graphs on any categorified quantum set. Moreover, since the Jones projection $e=\cF_2(\id_{B_1})$ generating $B_1$ is the neutral element for the Schur multiplication $-\star-$ in $\End_{A-A}(B)\cong A'\cap B_1,$ it is therefore justified to regard {\bf any categorified graph $\cG=({}_A{B}_A, \hat{T})$ as a sub-graph of the complete categorified graph} $({}_A{B}_A,e)$ which also contains the trivial categorified graph. Formally, this means $e\star \hat{T} = \hat{T} =\hat{T}\star e,$ or equivalently, the tautological $\id_{B_1}\circ T = T = T\circ \id_{B_1}.$

\begin{ex}[{\bf some $\Gamma$-equivariant  graphs}]\label{ex:GammaGraphs}
    We continue with the considerations from Example \ref{ex:GroupHyperfiniteCommutants}, where we used the Fourier transform to construct $\star$-idempotents on a finite group $\Gamma,$ or rather $\ell^\infty(\Gamma)$, corresponding to $\circ$-idempotents on $\bbC[\Gamma].$ 
    We shall now discuss cases where one can write these $\star$-idempotents explicitly. 

    The (abstract) semisimple decomposition of $\bbC[\Gamma]$ as a finite dimensional C*-algebra corresponds to the decomposition of the left regular representation of $\Gamma$ into its irreducible representations: 
    $$\ell^2(\Gamma)\cong \bigoplus_{\pi\in\mathsf{Irrep}(\Gamma)}\pi^{\oplus d_\pi},\quad \text{ with }\quad |\Gamma|= \sum_{\pi\in\mathsf{Irrep}(\Gamma)} d_\pi^2\quad \text{ and }\quad  \bbC[\Gamma]\cong \bigoplus_{\pi\in\mathsf{Irrep}(\Gamma)} \mathsf{Mat}_{d_\pi}(\bbC).$$
    
    With this decomposition in hand, we know everything about the group algebra $\bbC[\Gamma]$, which includes knowledge of all its $\circ$-idempotents (modulo similarity) and via the Fourier transform, we can obtain {\bf all adjacency matrices supported on the quantum set $\cR\subset \cR\rtimes\Gamma.$ }
    A problem with this approach is that obtaining a concrete isomorphism $M'\cap M_2\cong \bbC[\Gamma] \cong \bigoplus_{\pi\in\mathsf{Irrep}(\Gamma)} \mathsf{Mat}_{d_\pi}(\bbC)$ requires full knowledge of the conjugacy classes of $\Gamma,$ as the indicator functions of conjugacy classes yield the central projections in $\bbC[\Gamma]$. 
    Alternatively, one could cook up such an isomorphism with full knowledge of the representation theory of $\Gamma$ together with a fiber functor; this is  concrete Hilbert spaces and equivariant maps between them, as well as chosen orthonormal bases for each. This information allows one to realize $\bbC[\Gamma]$ acting on $\ell^2(\Gamma),$ which is  generated by the matrix coefficients of irreducible representations. 
    Either approach forces us to make unnatural choices and/or perform computations which will heavily depend on the data of $\Gamma.$
        
    In the finite abelian case, the Fourier transform gives an algebra isomorphism $\bbC[\Gamma]\cong \ell^\infty(\widehat{\Gamma}),$ and so all adjacency matrices on $\cR\subset\cR\rtimes\Gamma$ correspond to subsets of $\widehat\Gamma.$  In a nutshell, one just has to understand the abelian dual group $\hat \Gamma$. 
\end{ex}

We now recover Cayley graphs from this formalism. 

\begin{ex}[{\bf Concrete Schur idempotents in groups and their Cayley graphs}]\label{ex:GroupsAndCayleyGraphs}
    We will continue using the notation from Example \ref{ex:GroupHyperfiniteCommutants}.
    In the spirit of the concreteness, we shall show how to use our constructions to yield every Cayley graph $\mathsf{Cay}(\Gamma, S)$ explicitly, with $\emptyset\neq S\subseteq \Gamma.$ 
    The trick is to focus on the concrete categorified set ${_M}{M_1}_M$ instead of ${_N}M_N$ given by the Jones basic construction (Definition \ref{defn:JonesProjBasicConst}) corresponding to the dual subfactor $M\subset M_1$.
    
    A quantum graph $\cG$ on ${_M}{M_1}_M$, is given by a Schur idempotent $T\in \End({_M}{M_1}_M)\cong M'\cap M_2$. 
    We know from Example \ref{ex:GroupHyperfiniteCommutants}, that $$(M'\cap M_2, \star)\cong(N'\cap M_1, \circ)\cong \ell^\infty\Gamma.$$ 
    Thus, for some $\circ$-idempotent $\check T\in N'\cap M_1, $ we have $T = \cF_2^{-1}(\check{T}).$ 
    Notice that $\check T$ corresponds to some nonempty subset $S\subseteq \Gamma$, and so quantum graphs on the quantum set ${_M}{M_1}{_M}$ are parameterized by $2^\Gamma\setminus\{\emptyset\}.$
    Thus, we write $\check T=1_S = \sum_{s\in S}P_s$ as a sum of indicators of point-masses. 
    (Compare with \cite[Theorem D]{2023arXiv230615315W}, where $\Gamma$ is a discrete quantum group, with the difference that we do not demand $S$ be generating.)

    We shall now compute the inverse Fourier transform of each $P_s,$ for which we will use the usual calculus for the 2-category generated by the bimodule $X={_M}M_N$, with conjugate $N$-$M$ bimodule $\overline{X}$. We have that 
    \begin{align}
        \begin{aligned}
            \coev{X}: M&\to X\otimes\overline{X}\\
            m&\mapsto \sum_{g\in \Gamma} u_g\otimes \overline{u_g}\ m
        \end{aligned}
        \quad\text{with adjoint}\quad
        \begin{aligned}
            \coev{X}^\dag:X\otimes \overline{X}&\to M\\
            \eta\otimes\overline{\xi}&\mapsto {_M}\langle \eta,\ \xi\rangle^X = \eta\xi^*.
        \end{aligned}
    \end{align}
    Here, to ease notation we have returned to the convention that the round tensor $\otimes$ is balanced over $N$, while $\boxtimes$ means balancing over $M$. 
    We notice that, the following complimentary solution to the conjugate equations and its adjoint given by 
    $$\ev{X}(\overline{\eta}\boxtimes\xi)=\rip{\eta}{\xi}{N}\hspace{.5cm} \text{ and }\hspace{.5cm} \ev{X}^\dag(n)= \overline{u_e}\boxtimes u_e\ n$$ 
    might not yield standard/balanced solutions. 
    This is easy to fix, but we shall not need to. In this way we do not have to carry the factors of $|\Gamma|^{\pm 1/2}.$ 
    We shall moreover use the unitor morphism and its adjoint  
    \begin{align*}
        \begin{aligned}
            \lambda_X:\overline{X}\boxtimes X&\to M\\
            \overline{m_1}\boxtimes m_2&\mapsto m_1^*m_2
        \end{aligned}
        \qquad \text { with adjoint }\qquad
        \begin{aligned}
            \lambda_X^\dag:M&\to \overline{X}\boxtimes X\\
            m&\mapsto \overline{u_e}\boxtimes m,
        \end{aligned}
    \end{align*}
    turning a single $M$ strand into a double $\overline{X}\boxtimes X$ strand and \emph{vice versa}.

    The inverse quantum Fourier transform then becomes
    \begin{align}
        \begin{aligned}
        \cF_2^{-1}(\lambda_X^\dag\circ P_s\circ \lambda_X): X\otimes\overline{X}&\to X\otimes\overline{X}\\
        u_a\otimes\overline{u_b} &\mapsto u_{as}\otimes\overline{u_{bs}}    
        \end{aligned}
        \qquad\qquad \forall s\in S, \forall a,b\in \Gamma.
    \end{align}
    Thus, as a matrix indexed over $\Gamma\times \Gamma,$ we have 
    $$\left[\check{P_s}\right]_{a\overline{b}}^{c\overline{d}}= \delta_{c=as}\cdot \delta_{d=bs},$$
    mapping the edge $u_a\otimes\overline{u_b}$ to the edge $u_{sa}\otimes\overline{u_{sb}}.$
    Notice that if we let $a=b$ and $c=d,$ we recover the adjacency matrix of  Cayley graph of $\Gamma$ generated by $\{s\}$ on the nose. 
    Since $1_S= \sum_{s\in S}P_s$, and the Fourier transform is a linear map, {\bf we concretely recover the adjacency matrix of the Cayley graph $\mathsf{Cay}(G,S).$} 
\end{ex}

In fact, Examples \ref{ex:GammaGraphs} and  \ref{ex:GroupsAndCayleyGraphs} is a manifestation of a much more general phenomena where  one or both of the relative commutants $A'\cap B_1$, $B'\cap B_2$ of a finite-index inclusion $A\subset B$ is an abelian C*-algebra. 
Many such examples can be constructed from Hadamard matrices/spin models and we will further analyze them in forthcoming work. 

\medskip

From Examples \ref{ex:GroupHyperfiniteCommutants},   \ref{ex:GammaGraphs} and \ref{ex:GroupsAndCayleyGraphs} it is already apparent that categorified graphs are equivariant with respect to some symmetry. 
In Example \ref{ex:FinQuGroupoids} for instance, this can be taken as far as finite quantum groupoids. 
More generally, modelling quantum graphs over inclusions $A\subset B$, realizing $B\cong A\rtimes Q$ for some Q-system $Q$ in the UTC generated by ${}_AB_A\in \fgpBim(A)$ (cf \cite{MR4419534}), one can therefore think of all categorified graphs as generalized Cayley graphs, equivariant with respect to a quantum symmetry. 

We shall now incorporate the finite dimensional quantum graphs studied in the literature into our framework.
\begin{ex}[\textbf{All simple finite directed quantum graphs on $(M_N(\bbC),\ \omega)$}]\label{ex:QuGraphsMN} 
Consider the inclusion $\bbC=A\subset M_N =B$ equipped with the conditional expectation $E=\omega,$ which is the faithful state determined by $\omega(\xi)=\Tr(W\xi),$ where $W=\mathsf{diag}\left((q_i)_1^N\right).$ Here, the weights $(q_i)_1^N$ assemble into a nondegenerate probability distribution on $\{1, 2...,N\}$; i.e. for each $i$ have $q_i\in(0,1)$ and $\sum_1^N q_i = 1.$ Notice $\omega$ is tracial if and only if it corresponds to the uniform probability distribution. 

Equipped with this data, consider the GNS construction $L^2(B,\omega)$ with right inner product $$\rip{\eta}{\xi}{A}:=\omega(\eta^*\xi).$$ This is canonically a right C* $A$-correspondence with a left $B$ action given by multiplication. We denote this correspondence by $X:={_B}X_A:={_B}L^2(B,\omega)_A,$ and we shall now describe $X$ in detail.  
A computation reveals that $$\left\{u_{ab} = q_b^{-1/2}e_{ab}\right\}_{a,b=1}^N$$ forms a right orthonormal Pimsner-Popa basis for $X_A.$ 
The Watatani index $\mathsf{Ind}_W(E)=[B:A]$ is given by 
$$[B:A]= \sum_{a,b=1}^Nu_{ab}u_{ab}^*=\left(\sum_{i=1}^Nq_i^{-1}\right)I_N\in [1,\infty)I_N\subset M_n(\bbC).$$
Notice that when $W$ is the uniform distribution, or equivalently, when $\omega=\tr$, then $[B:A]= N^2$ attains its minimal value.
A left Pimsner-Popa basis for ${_B}X$ is $\{I_N\}$ with left action by multiplication and left inner product 
$$\lip{B}{\eta}{\xi}:=[B:A]^{-1/2}\cdot \eta\xi^*.$$  
To ease notation, when we take a relative tensor product balanced by $\bbC$ we sometimes write $\otimes$ instead of $\boxtimes_{\bbC}$, and when it is balanced by $M_N,$ we write $\boxtimes$ instead of $\boxtimes_{M_N}$. 
Viewing $m,$ the multiplication map on $B\otimes B$ as an operator, $L^2(B,\omega)\otimes L^2(B,\omega)\to L^2(B,\omega)$ one can see from direct computation that its adjoint  $m^\dag:L^2(B,\omega)\to L^2(B,\omega)\otimes L^2(B,\omega)$ is determined by 
$$m^\dag(u_{ab})= \sum_{i=1}^Nq_i^{-1/2}u_{ai}\otimes u_{ib}.$$
Thus, the {\bf delta form condition} (ie $\delta^2\id_B = m\circ m^\dag$) on $(B,\omega)$ holds if and only if $\delta^2 = \Tr(W^{-1}) = \sum_{i=1}^Nq_i^{-1}=[B:A].$

Recall from Equation \ref{eqn:DualityMaps} that the duality maps for ${_B}X_A$ are given by 
\begin{align*}
    &\ev{X}: {_A}\overline{X}\boxtimes_B X_A\to {_A}A_A \qquad \qquad\qquad\text{ and } \qquad&& \coev{X}:{_B}B_B\to {_B}X\boxtimes_A \overline{X}_B\\
    &\hspace{2cm}\overline{\eta}\boxtimes\xi\mapsto \omega(\eta^*\xi)= \rip{\eta}{\xi}{A}, && \hspace{2cm}b\mapsto \left(\sum_{i,j=1}^N u_{ij}\boxtimes\overline{u_{ij}}\right)\lhd b,
\end{align*}
and whose adjoints are given by as shown by a direct computation: 
\begin{align}\label{eqn:evcoevAdjoints}
    &\ev{X}^\dag: {_A}A_A\to {_A}\overline{X}\boxtimes_B X_A \qquad \qquad\text{ and } && \coev{X}^\dag: X\boxtimes_A \overline{X}\to {_B}B_B\nonumber \\
    &\hspace{1.6cm} a\mapsto a\rhd [B:A]^{1/2}\sum_{i=1}^N q_i\ \overline{u_{ii}}\boxtimes u_{ii}, && \hspace{1.9cm} \eta\boxtimes\overline{\xi}\mapsto [B:A]^{-1/2}\eta\xi^* =\lip{B}{\eta}{\xi}.
\end{align}

Notice that $\ev{X}\circ\ev{X}^\dag(1)=[B:A]^{1/2}$ is the left dimension, and $\coev{X}^\dag\circ\coev{X}(1)=\sum_{i=1}^N1/q_i=[B:A]^{1/2}$ is the right dimension of $X$, and so our choices indeed give the unitary (balanced) dual $\overline{X}$ for $X$. 
Thus, the {\bf (left/right) quantum dimension} (see discussion before Example \ref{ex:Pants}) of ${_B}X_A$ is $$d_X=\sqrt{[B:A]^{1/2}[B:A]^{1/2}}= [B:A]^{1/2}=\delta\in[1,\infty),$$
illustrating the categorical significance of  the $\delta$-form condition.

We are now ready to describe the first few relative commutants of $\bbC\overset{\omega}{\subset}M_N(\bbC)$ (ie its standard invariant) using the alternating tensor powers of $X$ and $\overline{X}$. (Recall Construction \ref{construction:HigherRelComs}.) 
We shall then use this characterization to obtain all Schur idempotents on the quantum set $(B,\omega)$ through the quantum Fourier transform.  
We compute: 
\begin{equation}\label{ex:MNOmegaStdInv}
\begin{tikzcd}[column sep=-0.5em]
     &  \overbrace{\cong M_N}^{\overbrace{ \End\left({_\bbC}\overline{X}{_{M_N}}\right) }^{\cong A'\cap B }}
     \hookrightarrow
     & \overbrace{\cong \overline{M_N}\boxtimes I_N\otimes\overline{I_N}\boxtimes M_N}^{\overbrace{\End\left({_\bbC}\overline{X}\boxtimes X{_\bbC}\right)}^{\cong A'\cap {B_1}}} 
     & {}\\
     &  \underbrace{\cong\bbC}_{\underbrace{\End\left({_{M_N}}L^2(B,\omega){_{M_N}}\right)}_{\cong B'\cap B}}
     \arrow[hook]{u}
     \hookrightarrow
     &  \underbrace{\cong M_N}_{\underbrace{\End\left({_{M_N}}X_{_{\bbC}}\right)}_{\cong B'\cap {B_1}}}
     \arrow[hook]{u}
     \hookrightarrow
     & \underbrace{\cong I_N\otimes \overline{I_N}\boxtimes M_N \otimes \overline{M_N}}_{\underbrace{\End\left({_{M_N}}X\otimes \overline{X}{_{M_N}} \right).}_{\cong B'\cap {B_2}} }
     \arrow["\cF_2"']{ul}
\end{tikzcd}
\end{equation}

We now compute the quantum Fourier transform explicitly using that all $B$-$B$ bilinear endomorphisms of $X\boxtimes_A \overline X$ are $A$-compact---even finite rank---, and is thus  generated as the linear span of $\{|u_{ab}\rangle\langle u_{\alpha\beta}| \}_{a,b,\alpha,\beta=1}^N$ under the isomorphisms highlighted above. 
On a general endomorphism $T\in\End(X\otimes \overline{X})$, whose coordinatized expression is 
$$T(u_{\gamma\theta}\otimes\overline{u_{\zeta\sigma}}) = \sum_{\lambda,\rho,\mu,\nu=1}^N T^{\lambda\rho\ \mu\nu}_{\gamma\theta\ \zeta\sigma}\cdot u_{\lambda\rho}\otimes\overline{u_{\mu\nu}}.$$
With the input from the maps in \ref{eqn:evcoevAdjoints}, the quantum Fourier transform (Definition \ref{defn:FT}) is given by 
\begin{align}\label{eqn:CFTMNOmega}
    \cF_2 ( T ) (\overline{ u_{\alpha\beta}} \boxtimes  u_{\eta\kappa}) &=\delta\cdot\sum_{i,\lambda,\rho=1}^N q_i T^{\lambda\rho\ \eta\kappa}_{ii\ \alpha\beta}\cdot\overline{u_{ii}}\boxtimes u_{\lambda\rho}. 
\end{align}
In particular, if $T= (|u_{a'b'}\rangle\langle u_{c'd'}|)\otimes(|\overline{u_{ab}}\rangle\langle \overline{u_{cd}}|)$ we obtain 
\begin{align}\label{eqn:QFTMNCompacts}
    \cF_2 (|u_{a'b'}\rangle\langle u_{c'd'}|)\otimes(|\overline{u_{ab}}\rangle\langle \overline{u_{cd}}|) (\overline{u_{\alpha\beta}}\boxtimes u_{\eta\kappa}) = \delta\cdot\delta_{cd=\alpha\beta}\cdot\delta_{ab=\eta\kappa}\cdot\delta_{c'=d'}\cdot q_{c'}\cdot\overline{u_{c'c'}}\boxtimes u_{a'b'}.
\end{align}

Recall that the quantum Fourier transform $\cF_2$ maps a $\circ$-idempotent $T\in \End\left({_{M_N}}X\otimes \overline{X}_{{M_N}}\right)$ to a Schur idempotent $\hat T\in \End({_\bbC}\overline{X}\boxtimes_{M_N}X_{\bbC})$, which by Lemma \ref{lem:FTIso} is in fact a bijection. 
Consequently, Equation \eqref{eqn:CFTMNOmega} for $\cF_2$  yields all finite simple quantum graphs on $(M_N, \omega),$ when $T\in \End({}_{M_N}X\otimes \overline{X}_{M_N})\cong \overline{M_N}\otimes M_N$ is an idempotent. 
We shall now further exploit Equation \ref{eqn:QFTMNCompacts} to {\bf explicitly obtain all the quantum adjacency matrices} with respect to our chosen Pimsner-Popa basis $\{u_{\eta\kappa}\}_{\eta,\kappa=1}^N$. 
To do so, use the {\bf unitor} and its adjoint map, to transform a double strand labelled $\overline{X}\boxtimes X$ into a single strand labelled by $M_N$ and vice versa. Namely: 
\begin{equation*}
    \begin{aligned}
        \lambda_X:\overline{X}\boxtimes X&\to M_N\\
        \overline{b_1}\boxtimes b_2&\mapsto b_1^*b_2
    \end{aligned}
    \qquad\text{ and }\qquad
    \begin{aligned}
        \lambda_X^\dag:M_N&\to \overline{X}\boxtimes X\\
        b&\mapsto \overline{I_N}\boxtimes b
    \end{aligned}
\end{equation*}
These are easily seen to be unitary and mutual inverses. 
Composing Equation \eqref{eqn:QFTMNCompacts} with $\lambda_X\circ-\circ\lambda_X^\dag$ yields the $(\eta\kappa,\ \eta'\kappa')$-entry of the operator from Equation \eqref{eqn:QFTMNCompacts}:
\begin{equation}\label{eqn:AdjMatMN}            
    \left[\cF_2\Big((|u_{a'b'}\rangle \langle u_{c'd'}|)\otimes|\overline{u_{ab}}\rangle\langle\overline{u_{cd}}|\Big)\right]_{\eta\kappa}^{\eta'\kappa'}= \delta\cdot \delta_{ab=\eta\kappa}\cdot\delta_{c=d}\cdot\delta_{c'=d'}\cdot\delta_{c'=a'}\cdot\delta_{a'b'=\eta'\kappa'}\cdot q_c^{1/2}q{_a'}^{1/2}.
\end{equation}
This will be used in forthcoming examples.

We now verify directly that $\cF_2(\id_{X\otimes\overline{X}})= (\id_{X}\otimes\id_{\overline{X}})\hat{}$ is a scalar multiple of the Jones projection. 
On the one hand, the Jones projection $e\in\End(\overline{X}\boxtimes_{M_N}X)\cong A'\cap B$ coming from the projection on $L^2(B,\omega)$ onto $L^2(\bbC I_N,\omega)$ is given by:
$$e(\overline{u_{\alpha\beta}}\boxtimes u_{\eta\kappa}) = \delta^{-1}\ev{X}^*\circ\ev{X}(\overline{u_{\alpha\beta}}\boxtimes u_{\eta\kappa}) = \delta_{\alpha\beta=\eta\kappa}\sum_{\ell=1}^N q_{\ell}\ \overline{u_{\ell\ell}}\boxtimes u_{\ell\ell}.$$
On the other, in coordinates, we have that 
$(\id_{X\otimes\overline{X}})_{\gamma\theta\ \zeta\sigma}^{\lambda\rho\ \mu\nu}= \delta_{\gamma\theta=\lambda\rho}\cdot\delta_{\zeta\sigma=\mu\nu},$ 
so we obtain 
$$\cF_2(\id_{X}\otimes\id_{\overline{X}})(\overline{u_{\alpha\beta}}\boxtimes u_{\eta\kappa}) = \delta\cdot\sum_{i,\lambda,\rho=1}^N q_i (\delta_{ii=\lambda\rho}\cdot\delta_{\eta\kappa=\alpha\beta})\cdot\overline{u_{ii}}\boxtimes u_{\lambda\rho}= \delta\cdot\delta_{\eta\kappa=\alpha\beta}\sum_{i=1}^N q_i\overline{u_{ii}}\boxtimes u_{ii},$$
recovering the Jones projection up to the scalar $\delta.$
Similarly, if we start with the map $\coev{X}\circ\coev{X}^\dag\in\End(X\otimes\overline{X}),$ it is easy to see $\cF_2(\coev{X}\circ\coev{X}^\dag)=\id_{\overline{X}X}.$  Diagrammatically, this follows from the conjugate equations for $X.$ 
\end{ex}

We shall now use the techniques described above to recover the (non-)tracial undirected reflexive simple quantum graphs on $M_2$ classified in \cite{MR4481115, MR4514486}. 
It should be clear now that this method can take us beyond $\delta$-forms $\bbC\overset{\psi}{\subset} M_2(\bbC)$, whenever the computation of either relative commutant $B'\cap B_2$ or $A'\cap B_1$ is feasible. It is however important to demonstrate the use of it in the smallest quantum spaces, so we now unpack Example \ref{ex:QuGraphsMN}:
\begin{sub-ex}[{\bf Undirected reflexive (non-tracial) quantum graphs on $(M_2,\omega)$}]
    Here we will continue using the notation and conventions from Example \ref{ex:QuGraphsMN}. 
    When $N=2$, the standard invariant for the inclusion $\bbC\overset{\omega}{\subset}M_2$ from Equation \ref{ex:MNOmegaStdInv} indicates the quantum Fourier transform becomes a linear isomorphism 
    $$\cF_2:\underbrace{\End_{M_2-M_2}\left(M_2\otimes_{\bbC} \overline{M_2}\right)}_{\cong M_2 \otimes_\bbC M_2} 
    \longrightarrow
    \underbrace{\End_{\bbC-\bbC}\left(\overline{M_2}\boxtimes M_2\right)}_{\cong \overline{M_2} \otimes_\bbC M_2}
    $$
    whose action on idempotents on $M_2\otimes \overline{M_2}$ we now describe using Equation \ref{eqn:AdjMatMN}. For $a,b,a',b'\in\{1,2\}$ we let $E_{ab\ a'b'}$ be the $4\times 4$ matrix whose $ab\ a'b'$ entry equals $1$ and filled with zeros elsewhereso the $E_{ab\ a'b'}$ are a system of matrix units for $M_2\otimes M_2$.
    We then obtain that, if $c\in{1,2},$
    $$
            [\cF_2((|u_{a'b'}\rangle\langle u_{a'a'}|)\otimes|\overline{u_{ab}}\rangle\langle\overline{u_{cc}}|)]=\delta\cdot E_{ab\ a'b'} \sqrt{q_cq_{a'}}.
    $$
    Matsuda \cite[Theorems 3.1 and 3.5]{MR4481115} classified all undirected reflexive quantum graphs on $(M_2,\omega)$ and gave representatives of their quantum adjacency matrices modulo quantum (and classical) isomorphism. 
    In our format, these four matrices look like: 
    \begin{align}\label{eqn:quGraphsM2}
        &\begin{aligned}
            \hat{A}_1=\begin{pmatrix}
                1 &0 &0 & 0\\
                0 &1 &0 & 0\\
                0 &0 &1 & 0\\
                0 &0 &0 & 1\\
            \end{pmatrix},
        \end{aligned}
        \qquad
        &&\begin{aligned}
            \hat{A}_2=\begin{pmatrix}
                q_2^{-1} &0 &0 & 0\\
                0 &0 &0 & 0\\
                0 &0 &0 & 0\\
                0 &0 &0 & q_1^{-1}\\
            \end{pmatrix},
        \end{aligned}\\
        &\begin{aligned}
            \hat{A}_3&=\begin{pmatrix}
                1 &0 &0 & (q_1q_2)^{-1/2}\\
                0 &1 &0 & 0\\
                0 &0 &1 & 0\\
                (q_1q_2)^{-1/2} &0 &0 & 1\\
            \end{pmatrix},
        \end{aligned}
        \qquad
        &&\begin{aligned}
            \hat{A}_4=\begin{pmatrix}
                q_2^{-1} &0 &0 & (q_1q_2)^{-1/2}\\
                0 &0 &0 & 0\\
                0 &0 &0 & 0\\
                (q_1q_2)^{-1/2} &0 &0 & q_1^{-1}\\
            \end{pmatrix}.
        \end{aligned}\nonumber
    \end{align}
\end{sub-ex}
Here, $\hat{A}_1$ is the trivial graph and $\hat{A}_4$ is the complete graph. In Example \ref{ex:QuGraphsMN} we have already shown that $\hat{A}_1$ is a scalar multiple of $\cF_2(\delta^{-1}\coev{X}\circ\coev{X}^\dag)$, and that $\hat{A}_4=\delta^{-1}\cF_2(e).$ 
(Recall that for us, $\delta= (q_1)^{-1}+(q_2)^{-1} = 1/(q_1q_2)$, and should not be confused with Matsuda's $\delta':=q+q^{-1}=(q_1q_2)^{-1/2},$ where $q=q_2^{1/2}q_1^{-1/2}.$) 
A direct application of Equation \eqref{eqn:AdjMatMN} shows that we can obtain $\hat{A}_2$ from an idempotent as 
$$\hat{A}_2= \left[\dfrac{1}{\delta q_1q_2}\right]\cF_2\Big((|u_{11}\rangle\langle u_{11}|\otimes |\overline{u_{11}}\rangle\langle\overline{u_{11}}|) + (|u_{22}\rangle\langle u_{22}|\otimes |\overline{u_{22}}\rangle\langle\overline{u_{22}}|)\Big).$$ 
Finally, 
$$\hat{A}_3 = \cF_2\Big(\delta^{-1}\coev{X}\circ\coev{X}^\dag + |u_{11}\rangle\langle u_{11}|\otimes|\overline{u_{22}}\rangle\langle \overline{u_{22}}| + |u_{22}\rangle\langle u_{22}|\otimes|\overline{u_{11}}\rangle\langle \overline{u_{11}}|\Big),$$
which are readily seen to be idempotents with mutually orthogonal ranges, and so the sum inside $\cF_2$ defines a projection.

\begin{ex}\label{ex:BiProjectionsCombinatorially}
    In Example \ref{ex:BiprojSubalgs}, we mentioned that biprojections carry a graph-theoretic interpretation. We shall illustrate this with the inclusion $\bbC= A\subset B=\bbC^6$. Indeed, this inclusion contains the intermeddiate subalgebras $C= \bbC^2$ and $D = \bbC^3$, and so there are biprojections $c,d\in End(\bbC^6)$ whose ranges are these algebras. With respect to the canonical basis of $B$, the simplest matricial presentation for the biprojections $c,d$ and the Jones projection $e\in \End(\bbC^6)$ give the adjacency matrices  
    \begin{align*}
    &\begin{aligned}
            2c=\begin{pmatrix}
                1 &1 &0 &0 &0 &0\\
                1 &1 &0 &0 &0 &0\\
                0 &0 &1 &1 &0 &0\\
                0 &0 &1 &1 &0 &0\\
                0 &0 &0 &0 &1 &1\\
                0& 0 &0 &0 &1 &1
            \end{pmatrix},
        \end{aligned}
        \ 
        &&\begin{aligned}
            3d=\begin{pmatrix}
                1 &1 &1 &0 &0 &0\\
                1 &1 &1 &0 &0 &0\\
                1 &1 &1 &0 &0 &0\\
                0 &0 &0 &1 &1 &1\\ 
                0 &0 &0 &1 &1 &1\\
                0 &0 &0 &1 &1 &1
            \end{pmatrix},
        \end{aligned}
        \ 
        &&&\begin{aligned}
            6e=\begin{pmatrix}
                1 &1 &1 &1 &1 &1\\
                1 &1 &1 &1 &1 &1\\
                1 &1 &1 &1 &1 &1\\
                1 &1 &1 &1 &1 &1\\
                1 &1 &1 &1 &1 &1\\
                1 &1 &1 &1 &1 &1
            \end{pmatrix},
        \end{aligned}
    \end{align*}  
    corresponding to 3 copies of the 2-clique $\bbK_2$, 2 copies of the 3-clique $\bbK_3$ contained in the 6-clique $\bbK_6 = (\bbC^6, e/6).$ 
\end{ex}

\begin{ex}[{\bf All categorified graphs on hyperfinite subfactors of index $< 4$}]\label{Ex:HypGraphsLeq4}
    Consider a finite index subfactor $A\overset{E}{\subset} B$, where $A\cong B\cong \cR$ are hyperfinite $\rm{II}_1$-subfactors and $E$ is the canonical trace preserving conditional expectation. 
    Hyperfinite subfactors of finite depth and index $[B:A]<4$ are determined by their principal graphs and fall under an $ADE$ classification \cite[\S4.6]{MR999799}. 
    (The principal graphs of a subfactor is a combinatorial invariant which describes the irreducible bimodules appearing as subobjects of repeated alternating tensor powers of a subfactors \emph{fundamental representation}, following the pattern of Construction \ref{construction:HigherRelComs}. See \cite{MR3166042} for a recent survey on the classification of subfactors.)
    Here for each $n\geq 3$ there is an $A_n$ subfactor of index $4\cos^2\pi/(n+1),$ a $D_{2n}$ subfactor of index $4\cos^2(\pi/(4n-2))$, and pairs of complex conjugate  subfactors of $E_6$ and $E_8$-type at index values $4\cos^2(\pi/12)\cong 3.7320...$ and $4\cos^2(\pi/30)= 3.9562...$, respectively. (No $D_{\sf{odd}}$ or $E_7$ subfactor exists \cite[Theorem 1.4.3]{MR999799}).

    The simplest nontrivial UTCs associated to the $A_n$-series are the {\bf Temperley-Lieb-Jones} $\TLJ(\delta)$ categories, with $\delta = 2\cos(\pi/n).$ 
    These can equivalently be presented as categories of pair noncrossing partitions with loop parameter $\delta$.
    There are finitely many isomorphism classes of simple objects known as the {\bf Jones-Wenzl projections}, labeled $1=f^{(0)}, f^{(1)},..., f^{(n-1)}$. See \cite{2017arXiv171007362E}, and references therein. 
    It is known that the other $DE$-type categories can be explicitly obtained from the $A_n$ categories and its modules. 
    The $E_6$ case is explicitly described in \cite{MR1259914}, giving an explicit description of $A'\cap B_1\cong \bbC e \oplus \bbC (1-e)$. See also \cite[\S5]{Bisch97} and \cite[\S4]{MR1470857}.

    A hyperfinite subfactor $A\subset B$ with $[B:A]<4$ is therefore a concrete realization of the depth-2 Q-system ${}_AB\boxtimes_B B_A\cong {_A}B_A\in \fgpBim(A),$ which decomposes as the direct sum of all the even vertices of the classical graph at depth $0$ and $2$, and thus automatically yields a finite quantum set. 
    Therefore, from the principal graphs we get $A'\cap B_1\cong B'\cap B_2 \cong \bbC^2$ as linear spaces; ie the 2-box spaces are linearly spanned by the identity and the Jones projection, which therefore are the only Schur idempotents (i.e. adjacency matrices) on $\End({}_AB_A)$. 

    Thus, the only possible connected simple undirected graphs on finite-depth hyperfinite subfactors of index $< 4$ are the complete and the trivial graphs.
\end{ex}

Hyperfinite subfactors have been completely classified up to index $5+1/4$ \cite{MR4565376}, and there are many other familiar constructions of explicit examples above this range, where there is little hope of a meaningful classification. 
What we attempt to clarify here is that any concrete example of a subfactor---hyperfinite or not---, where the 2-box spaces are computable, affords new examples of finite concrete categorified graphs.
For instance, it is known that there's an $A_\infty$ (non-hyperfinite necessarily) subfactor at any index value above $4$, which yields examples of complete categorified graphs on finite sets of arbitrarily large dimension whose only possible simple graphs are either complete or trivial. 
These $A_\infty$ subfactors are associated to the diagrammatic $\TLJ$ categories/representations of $SU_q(2),$ where moreover all Q-systems are known. 
There are 2-categorical generalizations of $\TLJ$ which are constructed from classical graphs, all of whose fiber functors are known (\cite{MR3420332, MR4098904}).

\medskip

We now show an explicit small example of a genuine categorified graph supported on the Ising category:
\begin{sub-ex}
   An $A_3$ subfactor is supported on the UTC Ising $\cC=\mathsf{Ising}$ with $\Irr(\mathsf{Ising})= \{\mathbf{1}, \tau, \sigma\}$ with fusion rules $\tau\otimes\tau = \mathbf 1, \tau\otimes \sigma = \sigma = \sigma\otimes \tau$ and $\sigma\otimes \sigma = \mathbf{1}\oplus\tau.$ 
   Moreover, we know that the connected Q-systems in 
   $\mathsf{Ising} $ are $\{\mathbf{1},\ Q\}$. 
   The decomposition of $Q$ is given by $Q\cong \sigma\otimes \sigma \cong \mathbf 1\oplus \tau,$ with quantum dimension $d_\sigma=\sqrt 2= 2\cos(\pi/4)$ corresponding to a subfactor of index $2$. \cite[Example 3.14]{MR3308880} 
   From Example \ref{Ex:HypGraphsLeq4}, it follows that the only possible graphs supported on the $A_3$ subfactor are (direct sums of) either the trivial or the complete graph. 
\end{sub-ex}

\section{A Quantum Variant of Frucht's Theorem}

A renowned result due to Frucht \cite{MR1557026} states that every finite group is exactly the automorphisms of some finite simple classical graph, and a natural question to ask is if every finite quantum group arises as the quantum automorphisms of some  graph.  The answer turns out to be no, as was shown by Banica and McCarthy \cite{MR4462380}.  Among other examples, the authors show that the quantum group $\hat{S_3}$ dual to $S_3$ is not the quantum automorphism group of any such graph.

The goal of this section is to highlight how our enriched context for graphs leads to a setting in which a variant of  Frucht's theorem holds for quantum group(oid)s. 
We shall now profit further from the abundance of examples of categorified graphs from Section \ref{sec:IntroCatQuSetGraph} and the generality  achieved so far by showing that {\bf every finite quantum groupoid arises as certain `\emph{quantum automorphisms}' of some categorified graph}.
To achieve this, in Theorem \ref{thm:QuFrucht}, we will exploit the connections with the theory of subfactors and its reconstruction tools. 
However, we will first need to establish what we mean by a quantum automorphism of a categorified graph in our setting and a few more preliminaries

We shall first outline aspects of the work of Nikshych and Vainerman \cite{MR1897159, MR1745634, MR1800792, MR1913440}, where it was established that finite-index $\rm{II}_1$-subfactors $N\subset M$ of finite depth (not necessarily irreducible: $N'\cap M=\bbC1$) correspond to weak Hopf C*-algebras and their crossed products by actions on such factors. 
In the following example, we extend the scope of the duality observed in Example \ref{ex:GroupHyperfiniteCommutants} in the much broader context of finite-depth subfactors and associated Hopf algebras.  For the sake of space, we omit almost all the definitions pertaining to finite quantum groupoids, and refer the reader to e.g., \cite{MR1745634, MR1913440} for details.
\begin{ex}[{\bf Finite depth \& finite index subfactors and finite quantum groupoids}]\label{ex:FinQuGroupoids}
    Let $$\bbH = (H,m, 1, \Delta, \varepsilon, S)$$ be a \emph{finite quantum groupoid} (ie a weak Hopf C*-algebra).
    This is, $H$ is a finite-dimensional C*-algebra with comultiplication, counit and antipode maps
    \begin{align*}
        &\Delta:H \to H\otimes H    && \varepsilon: H\to \bbC &&& S:H\to H.
    \end{align*}
    making $(H,\Delta, \epsilon)$ a $\ast$-coalgebra.   Here, $\Delta$ is a $*$-homomorphism that need not be unital, and $S$ is an anti-(co)homomorphism satisfying the compatibility conditions from \cite[\S2]{MR1800792}. 
    If $S$ satisfies $S^2=\id_H$ we say $\bbH$ is a weak Kac algebra.  Note that a finite quantum groupoid is a quantum group if and only if $\Delta $ is unital if and only if $\epsilon $ is an algebra homomorphism.

    The vector space $\hat{H}$ dual to $H$ admits the structure of a finite quantum groupoid denoted by 
    $$
    \hat{\bbH} = (\hat{H}, m_{\hat{H}} = \Delta^\dag, \varepsilon, \Delta_{\hat{H}}= m^\dag, 1, S^\dag),
    $$
    where the $\dag$ denotes the adjunction given by the pairing between $H$ and its dual. Then $\hat{\bbH}$ is the quantum groupoid dual to $\bbH.$
    
    Consider the following canonical commuting square, 
    \begin{equation}\label{eq:GpoidCommutingSq}
    \begin{tikzcd}[column sep=0em, row sep=-.1em]
     \hat{H} & \subset   & \hat{H}\rtimes H\\
    \cup & {} & \cup \\
     \hat{H}\cap H & \subset &  H 
     \end{tikzcd},
    \end{equation}
    where the crossed product $\hat{H}\rtimes H$ C$^\ast$-algebra is constructed as in \cite[Secion 8]{MR1913440}.
    Assuming the inclusions above are connected, by carrying out the Jones tower construction \ref{construction:HigherRelComs}, we obtain a finite index hyperfinite subfactor of finite-depth $N\subset M$, where one can take (\ref{eq:GpoidCommutingSq}) as the canonical commuting square of the inclusion $M\subset M_1$ (c.f. \cite[\S5]{MR1897159} for the Kac algebra case).
    Without loss of generality, one can assume depth-2 by repeating the basic construction sufficiently and identifying the original inclusion as an intermediate subfactor \cite[Corollary 4.2]{MR1800792}.
    Moreover, \cite[Theorem 5.8]{MR1800792} guarantees that 
    $$
    \cC_{N\subset M}\cong \Rep(\hat\bbH) 
    $$ 
    as unitary fusion categories, where the former is the tensor category generated by the finitely generated projective $N$-$N$ bimodules generated by $N\subset M$ and the latter that of finite dimensional modules over $\hat H$. 
    Similarly, we have 
    $$
    \cC_{M\subset M_1}\cong \Rep(\bbH). 
    $$
    In particular, the commuting square: 
    \begin{equation}\label{eq:CommSq}
    \begin{tikzcd}[column sep=0em, row sep=-.1em]
     N'\cap M_1 & \subset   & N'\cap M_2\\
    \cup & {} & \cup \\
     M'\cap M_1 & \subset &  M'\cap M_2 
     \end{tikzcd}.
    \end{equation}
    is isomorphic to (\ref{eq:GpoidCommutingSq}), where in fact, by \cite[\S3]{MR1745634}, the finite dimensional C*-algebras $N'\cap M_1$ and $M'\cap M_2$ have the structure of finite quantum groupoids and are dual to each other, with 
    $$
    \hat{H}\cong N'\cap M_1\qquad \text{ and }\qquad H\cong M'\cap M_2.
    $$ 
    This generalizes the Pontryagin duality we saw on Example \ref{ex:GroupHyperfiniteCommutants}.

    Furthermore, there is an action $H\curvearrowright M_1$ coming from the unique minimal action of $H$ on the hyperfinite $\rm{II}_1$-factor given by iterating the basic construction on the square (\ref{eq:GpoidCommutingSq}) (c.f. \cite[Corollary 5.6]{MR1897159} for the  weak Kac algebra case).
    We can thereafter re-express the algebras in the Jones Tower as
    $$
    \underbrace{M_1^H\cong M}_{\text{ fixed-points }}, \qquad 
     \underbrace{M_2\cong M_1\rtimes H}_{ \text{crossed product } }\qquad \text{ and } \qquad \underbrace{M_3\cong M_2\rtimes \hat{H}}_{ \text{crossed product } },
    $$
    and so on.
    Moreover, the duality maps between $\bbH$ and $\hat\bbH$ can be described explicitly in terms of the constituents of the tower.
\end{ex}

From the above we see that the Jones Tower is constructed by taking crossed products by alternating actions of $\bbH$ and $\hat{\bbH}$ over $M_1$. 
In particular, $\hat{\bbH}$ acts as the quantum automorphisms of the quantum set given by the crossed product $M_1\rtimes H\cong M_2.$ 
It is therefore meaningful to determine a precise sense in which $\hat{\bbH}$ acts as ``quantum \emph{graph} automorphisms'' on $\bbK_{M_2},$ the complete graph over $M_2.$ To this end, we lay out a general categorical description.

Given a C* 2-category $\cC$, we consider {\bf the 2-category of $\cC$-equivariant graphs} $\Graph_\cC$ to be the 2-category whose objects are categorified graphs $\cG=(Q, \hat{T}),$ where $Q\in \cC(a\to a)$ is a Q-system in $\cC$, and $\hat{T}\in\End_{\cC}(Q)$ is a Schur idempotent. 
Next, we define the higher layers of $\Graph_\cC$ by generalizing Musto-Reutter-Verdon's 2-category of quantum graphs on finite-dimensional C*-algebras inspired on Gelfand's duality. 

Given $\cG=(Q,\hat{T}), \cG'=(Q', \hat{T'})\in \Graph_\cC,$ a 1-morphism is a \emph{quantum function} $(f,K),$ where $K\in \cC(a'\to a)$ and $f\in \cC(K\otimes_a Q\Rightarrow Q'\otimes_{a'}K)$ satisfies Relations (23) in \cite[Definition 3.11]{MR3849575} and also preserves the adjacency data as in Relations (51) in \cite[Definition 5.4]{MR3849575}. 
The compositional rule for quantum functions is also defined therein. 
In case $Q=Q'$ and $K=1_a$, $f$ is precisely a co-unital $*$-co-homomorphism on $Q$.
We say $(f,K)$ is a \emph{quantum isomorphism} if it also satisfies Relations (34) in \cite[Definition 4.3]{MR3849575}. 
When $K=1_a$ as before, a quantum isomorphism is additionally a unital $*$-endomorphism of $Q.$
The space of 2-morphisms between quantum functions $(f,K), (f',K')\in \Graph_\cC(\cG\to \cG')$ is given by those $h\in\cC(K\Rightarrow K')$ compatible with $f.$

We shall further explore the categorical structure and properties of $\Graph_\cC$ in a forthcoming paper. 
For now, we will state a useful lemma and the last condition we shall need to establish a categorified version of Frucht's theorem. 
\begin{lem}\label{lem:BraidingQuantumFunctions}
    Let $\cC$ be a unitary tensor category (viewed as a C* 2-category with one object) and $\cG = (Q, \hat T)\in \Graph_\cC$ be a graph. 
    Assuming further that $\cC$ has a unitary braiding $\beta=\{\beta_{a,b}:a\otimes b\to b\otimes a\}_{a,b\in\cC}$ (c.f. \cite[\S8.1]{MR3242743}), it follows that for all $K\in \cC$ we have 
    $$(\beta_{K,Q}, K)\in \Graph_{\cC}(\cG)$$
    is a quantum isomorphism of $\cG.$
    Moreover, the endomorphisms  of $(\beta_{K,Q}, K)$ are
    $$
    \Graph_\cC\left((\beta_{K,Q}, K)\Rightarrow (\beta_{K,Q}, K)\right) = \cC(K\to K).
    $$
\end{lem}
\begin{proof}
    Checking that $\beta_{K,Q}$ defines a quantum isomorphism follows from the naturality of $\beta$ in both entries as well as routine diagrammatic manipulations using the standard pictorical presentation of the braiding:  
    \begin{align*}
        \beta_{K,Q}=
    \tikzmath{
    \begin{scope}
        \clip[rounded corners=5pt] (-.5,-.5) rectangle (.5,.5);
        \fill[\AColor] (-.5,-.5) rectangle (.5,.5);
    \end{scope}
    \draw (-.3,.5) node[above]{$\scriptstyle Q$} --(-.3,.3) .. controls ++(270:.2cm) and ++ (90:.2cm) .. (.3,-.3) --(.3,-.5) node[below]{$\scriptstyle Q$};
    \begin{scope}
        \clip[rounded corners=5pt] (-.1,-.1) rectangle (.1,.5);
        \fill[\AColor] (-.1,-.1) rectangle (.1,.1);
    \end{scope}
    \draw (-.3,-.5) node[below]{$\scriptstyle K$}--(-.3,-.3) .. controls ++(90:.2cm) and ++ (270:.2cm) .. (.3,.3) --(.3,.5) node[above]{$\scriptstyle K$};
    },
    \end{align*}
    granted by its defining properties. 

    The last assertion is also immediate from the naturality. 
\end{proof}

Whenever $(\cC, \beta)$ is a braided unitary tensor category, and $\cG\in \Graph_{\cC}$, it then makes sense to restrict to the subcategory of {\bf crossing quantum isomorphisms of $\cG$},
$$
    \Graph_{\cC}^\beta(\cG)\subset \Graph_{\cC}(\cG\to \cG), 
$$
generated by quantum isomorphisms of $\cG$ the form $(\beta_{K, Q}, K),$ where $Q\in \cC(a\to a)$ is the underlying quantum set of $\cG$ and $K\in \cC(a\to a)$ is arbitrary.

We are now ready to state the following theorem: 
\begin{thm}[{\bf a categorified Frucht theorem} (Theorem \ref{thmalpha:Fruchts})]\label{thm:QuFrucht}
    Let $\bbH$ be a quasitriangular (equivalently, if $\Rep(\bbH)$ is braided cf \cite[\S5.2]{MR1913440}) finite quantum groupoid. 
    Then there is a hyperfinite finite-index subfactor $N\subset M$ of finite depth such that if $\cC=\cC_{M\subset M_1}$ is the standard invariant of the basic construction $M\subset M_1$, and $\cG = \cK_{M\subset M_1}$ is the complete graph on the Q-system ${}_M{M_1}_M$ then
    $$
        \Graph_{\cC}^{\beta}(\cG)\overset{\otimes}{\cong}\ \cC \ \overset{\otimes}{\cong} \Rep(\bbH).
    $$
\end{thm}
\begin{proof}
    Take the finite-index finite-depth subfactor $N\subset M$ constructed in Example \ref{ex:FinQuGroupoids}, and its basic construction $M\subset M_1$, and consider $\cC$ and $\cG\in \Graph_\cC$ as in the statement.
    (Note that we could as well take any $M$-$M$ sub-bimodule of ${}_M{M_1}_M$ generating $\cC$.)

    By \cite[Theorem 5.8]{MR1800792}, we already know that $\Rep(\bbH) \cong \cC$ as unitary fusion categories, so we need only to show the first equivalence. 

    Let $\beta$ be the braiding on $\Rep(\bbH)\cong \cC_{M\subset M_1}$ obtained from quasitriangularity. 
    Automatically, $\beta$ is unitary, as $\Rep(\bbH)$ is a unitary fusion category, by \cite{MR3239112}. 
    Consider the functor 
    \begin{align}
        \mathsf{B}:\cC &\to \Graph_{\cC}^{\beta}(\cG) \\    
        K&\mapsto (\beta_{K,M_1},\ K)\nonumber\\
        (h:K\to K')&\mapsto h.\nonumber
    \end{align}
    By the naturality of $\beta$, and Lemma \ref{lem:BraidingQuantumFunctions}, it follows that $\mathsf{B}$ is well-defined; this is, $\beta$ is a quantum graph automorphism, and $h$ is a quantum isomorphism intertwiner. 
    This essentially shows that every $K\in\cC_{M\subset M_1}$ does give rise to some quantum bijection on $M_1.$
    
    Moreover, there is a forgetful functor
    \begin{align}
        \For: \Graph_{\cC}^{\beta}(\cG)&\to \cC\\
        (\beta_{K,M_1},\ K)&\mapsto K\nonumber\\
        (h: K\to K')&\mapsto h.\nonumber
    \end{align}
    It is immediate that $\For$ is faithful (i.e. injective at the level of morphisms). 
    Furthermore, by the naturality of $\beta,$ and since we have restricted to the appropriate category we see immediately that $\For$ is surjective on morphisms:
    $$
    \For[\Graph_{\cC}^{\beta}(\cG)(f\to f')]= \cC(K\to K').$$ 
    That is, $\For$ is full.
    Finally, since $\mathsf{B}$ is defined on all of $\cC,$ it follows that for every $K\in \cC,$ we have $\For(\beta_{K, M_1},K)= K,$ and so $\For$ is surjective. 
    We conclude by general abstract nonsense, as a fully-faithful (essentially) surjective functor gives an equivalence. 
    The monoidality-preserving is obvious. 
\end{proof}

We remark that there are plenty of examples of quasitriangular quantum groupoids, such as Drinfeld doubles of ordinary finite quantum groupoids (see \cite[\S5.3]{MR1913440}), modular tensor categories, or even group Hopf algebras $\bbC[\Gamma]$. 
Moreover, our  Theorem \ref{thm:QuFrucht} could be compared with \cite[Proposition 5.19]{MR3849575}, in case $\cC=\fdHilb.$

The assumption of quasitriangularity in our quantum Frucht's theorem might be non-essential.
We only used it to have a source of nontrivial quantum bijections on quantum sets, which would allow us to remember all the data of $\Rep(\bbH)$ after applying the Forgetful functor. 
At present, we do not know of a better way to turn these obstacles around, while perhaps there might be more general versions of Frucht's Theorem.

\begin{remark}
    We conclude this paper by drawing  some parallels between the classical Frucht theorem and its proof with our Theorem \ref{thm:QuFrucht}. 
    
    Let $\Gamma$ be a finite classical group, and a $S\subset \Gamma$ a generating symmetric subset. 
    First, consider the Cayley graph $\cG':=\mathsf{Cay}(\Gamma,S),$ regarded as a colored directed graph. 
    That is, numbering $S=\{s_1,..., s_N\}$, there is an $i$-colored edge from $g$ to $h$ if $gs_i=h.$ One can then show that $\Aut^{\textsf{col}}(\cG')$, the group of color-preserving automorphisms of $\cG'$ equals $\Gamma,$ by observing these correspond to multiplication by group elements. 

    The second step is to construct an undirected simple graph $\cG$ with $\Aut(\cG) \cong\Aut^{\textsf{col}}(\cG')\cong \Gamma$.
    To achieve this, for all $i=1,2,...,N$, replace each $i$-colored directed edge $s_ig\to h$, by a \emph{spike-like path} undirected graph consisting of a line of $i+2$ new vertices, where the first is connected to $g$, and the last to $h$ and all  intermediate new vertices are also connected to a degree-one additional vertex, except the second to last, which is connected to a path of length two as follows: 
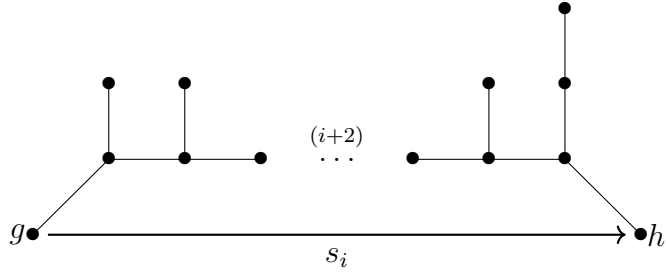
\begin{align*}
    \tikzmath{
    \draw (-4.2,0)node{$g$};
     \draw (-4,0)node{$\bullet$} (-4,0)--(-3,1) node{$\bullet$};
     \draw (-3,1) -- (-3,2) node{$\bullet$} (-3,1) -- (-2,1) node{$\bullet$};
     \draw (-2,1) -- (-2,2) node{$\bullet$} (-2,1) -- (-1,1) node{$\bullet$};
    \draw (0,1)node{$\hdots$};
    \draw (0,1.3)node{$\scriptstyle (i+2)$};
    \draw[thick, ->] (-3.8, 0) -- (3.8,0);
    \draw (0,-.3)node{$s_i$};
    \draw (4,0)node{$\bullet$} (4,0)--(3,1) node{$\bullet$};
     \draw (3,1) -- (3,2) node{$\bullet$} (3,1) -- (2,1) node{$\bullet$} ;
     \draw (2,1) -- (2,2) node{$\bullet$} (2,1) -- (1,1) node{$\bullet$} (3,2)--(3,3) node{$\bullet$};
     \draw (4.2,0)node{$h$};
     }
\end{align*}

    We now recall Example \ref{ex:GroupHyperfiniteCommutants} in light of Theorem \ref{thm:QuFrucht} and its proof.
    The choice of a symmetric generating set $S\subseteq \Gamma$ is manifested in our proof as follows: represent $S$ in $\ell^\infty(\Gamma)$ by the projection $1_S = \sum_{s\in S}\delta_s$. 
    Then using the quantum Fourier transform to yield adjacency data from $S$, the partition 
    $$
    \cF^{-1}(1_S) = \sum_{s\in S}\cF^{-1}(\delta_s) =:\sum_{s\in S}P_s
    $$
    yields a $\star$-idempotent in $\End_{M-M}(M_1)$ that keeps track of the coloring of edges by elements of $S.$ 
    That $S$ generates $\Gamma$ translates to 
    $$
    C^*(\{\cF^{-1}(\delta_s)= P_s\}_{s\in S}) \cong \End_{M-M}(M_1)
    $$
    with composition, and so the quantum set of vertices is the corresponding object ${}_M{M_1}_M.$ 
    This, in turn, ensures that the UTC generated by the vertices is $\Rep(\bbC[\Gamma])$. 
    The statement of Theorem \ref{thm:QuFrucht}---as well as the isomorphisms of the commuting squares (\ref{eq:GpoidCommutingSq}) and (\ref{eq:CommSq})---precisely tells us that the crossing quantum bijections of ${}_M{M_1}_M$ (automatically preserving the adjacency) are precisely the Hopf algebra $\bbC[\Gamma]$ we started with.
    And moreover that the $M$-$M$ bilinearity of these maps exactly captures that these automorphisms preserve the coloring by elements of $S$. 
    
    Even more directly, recalling that $\cF^{-1}(\delta_s)$ is the adjacency matrix of $\mathsf{Cay}(\Gamma, \{s\}),$
    that $C^*(\{\cF^{-1}(\delta_s)= P_s\}_{s\in S}) \cong \End_{M-M}(M_1)$ says that the naturally allowed endomorphisms of $M_1$ are traced back through the Fourier transform from the elements $s$ of the generating set $S.$
    This directly categorifies the first step in Frucht's proof in terms of the Pontryagin duality
     
    It is not clear to us how to achieve an analogue of the second step in the proof of Frucht's classical theorem that carefully enlarges the graph to blend the coloring into the adjacency data. However, it is also not clear if this is  necessary, or even meaningful, as the coloring and the adjacency simultaneously exist in dual layers in a C* 2-category of bimodules which are remembered by a finite quantum groupoid and its dual.
\end{remark}

\bibliographystyle{amsalpha}
\bibliography{bibliography}

\end{document}